\documentclass[a4paper,11pt,reqno]{amsart}
\usepackage{amssymb,amsthm,amsmath}
\usepackage{mathtools}
\usepackage{amsfonts}
\usepackage{amscd}
\usepackage{amssymb}
\usepackage{enumerate}
\usepackage{bbm}
\usepackage{graphicx}
\allowdisplaybreaks
\usepackage{color}
\usepackage{amsbsy}
\usepackage{graphicx}
\usepackage{amsthm}
\usepackage{amsmath}
\usepackage{amsxtra}
\usepackage{mathrsfs}
\usepackage{bbm}
\usepackage{dsfont}
\usepackage{enumitem}
\usepackage{comment}
\usepackage{xcolor}
\usepackage{float} % To use [H]
\usepackage{url}
\usepackage{soul}
\usepackage[hidelinks]{hyperref}

\usepackage{ifthen}
\hypersetup{
    colorlinks=true,
    linkcolor=red,      % for \ref and \eqref
    citecolor=red,       % for \cite
    urlcolor=cyan        % for \url
}

\newtheorem{theorem}{Theorem}[section]
\newtheorem{defi}[theorem]{Definition}

\newtheorem{corollary}[theorem]{Corollary}
\newtheorem{prop}[theorem]{Proposition}
\newtheorem{lemma}[theorem]{Lemma}

\numberwithin{equation}{section}
\definecolor{red}{rgb}{1.0, 0.0, 0.0}
\newcommand{\Bea}{\begin{eqnarray*}}
    \newcommand{\Eea}{\end{eqnarray*}}
\newcommand{\Be} {\begin{equation*}}
    \newcommand{\Ee} {\end{equation*}}
\newcommand{\be} {\begin{equation}}
    \newcommand{\ee} {\end{equation}}
\newcommand{\bea} {\begin{eqnarray}}
    \newcommand{\eea} {\end{eqnarray}}

\newcommand{\F}{\mathcal{F}}
\newcommand{\Hg}{\mathbb{H}}
\newcommand{\R}{\mathbb{R}}

\newcommand{\al}{\alpha}

\newcommand{\la}{\lambda}
\newcommand{\C}{{\mathbb C}}

\renewcommand{\Im}{\operatorname{Im}}

\usepackage{color}

\title[Uncertainty Principles on the Heisenberg Group]{Uncertainty Principles for the Strichartz Fourier transform on the Heisenberg Group}
\author[A. Dabra]{Arvish Dabra}
\address{
    Arvish Dabra:
    \endgraf
    Department of Mathematics
    \endgraf
    Indian Institute of Technology Delhi
    \endgraf
    Hauz Khas, New Delhi - 110016
    \endgraf
    India
    \endgraf
    {\it E-mail address:} {\rm arvishdabra3@gmail.com}
}

\author[A. Dasgupta]{Aparajita Dasgupta}
\address{
    Aparajita Dasgupta:
    \endgraf
    Department of Mathematics
    \endgraf
    Indian Institute of Technology Delhi
    \endgraf
    Hauz Khas, New Delhi - 110016
    \endgraf
    India
    \endgraf
    {\it E-mail address:} {\rm adasgupta@maths.iitd.ac.in}
}

\author[P. Gulia]{Prerna Gulia}
\address{
    Prerna Gulia:
    \endgraf
    Department of Mathematics
    \endgraf
    Indian Institute of Technology Delhi
    \endgraf
    Hauz Khas, New Delhi - 110016
    \endgraf
    India
    \endgraf
    {\it E-mail address:} {\rm prernagulia64@gmail.com}
}


\vspace{1cm}

\allowdisplaybreaks

\begin{document}

\begin{abstract}
In this article, we establish several fundamental uncertainty principles for the Strichartz Fourier transform on the Heisenberg group, including Benedicks’ theorem, the Donoho-Stark principle, the local uncertainty principle of Price, and a weak form of Beurling’s theorem. The Strichartz Fourier transform, introduced by Thangavelu (2023), provides a scalar-valued analogue of the classical operator-valued Fourier transform on the Heisenberg group. We first prove an analogue of Benedicks’ theorem asserting that a nonzero function and its Strichartz Fourier transform cannot both be supported on sets of finite measure. As a consequence, we obtain Nazarov’s uncertainty inequality. We then establish the Donoho-Stark principle, providing quantitative bounds on simultaneous concentration in space and frequency, and extend the local uncertainty principle of Price to this framework. Finally, we present a weak form of Beurling’s theorem for radial functions on the Heisenberg group.
\end{abstract}

\maketitle

\section{Introduction}

Uncertainty principles in harmonic analysis quantify the inherent limitations on the simultaneous localization of a function and its Fourier transform. Numerous formulations of these principles have been proposed in the literature; a comprehensive overview is provided in the survey by Folland and Sitaram \cite{FSsurvey}. We provide a concise overview of several uncertainty principles, including their proofs within our framework, their historical development and related contributions across various mathematical settings. Given the extensive literature in this field, it is impossible to cite all relevant works; we therefore highlight those contributions that have most profoundly influenced the development of this research area.

\subsection{Benedicks' theorem}
 In 1985, Benedicks \cite{benedicks-original-paper} proved that for any function \(f \in L^2(\mathbb{R}^n)\), if both the sets \(\{x\in \mathbb{R}^n:f(x) \neq 0\}\) and \(\{\xi\in \mathbb{R}^n:\widehat{f}(\xi)\neq 0\}\) have finite measure, then \(f\) must be identically zero. This result has inspired analogous formulations for various transforms in harmonic analysis. For a unimodular group \(G\) with unitary dual \(\widehat{G}\), Arnal and Ludwig \cite{arnal-ludwig} established the ``qualitative uncertainty principle" which states the following: if \(m\) and \(\widehat{m}\) denote the Haar measure on $G$ and Plancherel measure on \(\widehat{G}\), respectively and for \(f \in L^2(G)\), \(m\left(\{x\in G:f(x) \neq 0\}\right)<m(G)\) and 
$$
\int_{\widehat{G}}\operatorname{rank}\big(\widehat{f}(\lambda)\big) \, d\widehat{m}(\lambda)<\infty,
$$
then \(f=0\). When \(G\) is identified with the Heisenberg group $\mathbb{H}^n$, the aforementioned conditions imply that the Fourier transform of the function is supported on a set of finite measure and that \(\widehat{f}(\lambda)\) is a finite-rank operator. In 2010, Narayanan and Ratnakumar \cite{Narayanan-Ratnakumar} established that if \(f \in L^1(\mathbb{H}^n)\) is supported in \(B \times \mathbb{R}\), where \(B\) is a compact subset of \(\mathbb{C}^n\) and \(\widehat{f}(\lambda)\) is a finite-rank operator, then \(f=0\).
For further developments and related results on the Heisenberg group, see \cite{Price-Sitaram-transforms_on_LCT_groups,Sitaram-Sundari-Thangavelu-uncertainty_principles}.

Furthermore, in 2017, Vemuri \cite{Vemuri-benedicks} established Benedicks' theorem for the Weyl transform. More recently, in 2022, Ghosh and Srivastava \cite{somnath-benedicks_motion_group} proved a version of the theorem for the Heisenberg motion group and the Quaternion Heisenberg group. 

One can also refer to the works of Ghobber and Jaming \cite{jaming-Ghobber-bessel_transform,jaming-integral_operators} on several uncertainty principles for the Fourier-Bessel transform and integral operators, which have motivated various subsequent studies, see \cite{Daher-Opdam_Cherednik_transform,Negzaoui-Laguerre_transform_benedicks}.

\subsection{Donoho-Stark principle}
In the work of Donoho and Stark \cite{Donoho-stark_paper}, it was shown that in \(\mathbb{R}^n\), if a function \(f\) and its Fourier transform \(\widehat{f}\) both have unit norm and ``practically vanish" outside the sets \(T\) and \(W\), then there exists a lower bound on the product of the measures of these sets. The exact definition of ``practically vanishing" and the expression of the lower bound is given in terms of the concentration of the functions \(f\) and \(\widehat{f}\). Namely, for any \(\epsilon_T >0\), a function \(f\) is said to be \(\epsilon_T\)-concentrated on \(T\) if there exists a function \(g\) that vanishes outside \(T\) and satisfies \(\lVert f-g \rVert_2 \leq \epsilon_T\). The precise statement of the Donoho-Stark principle in $\mathbb{R}^n$ \cite{Donoho-stark_paper} is as follows:
\begin{theorem}
    Let $T$ and $W$ be measurable sets and suppose that the function $f$ and its Fourier transform $\widehat{f}$ are of unit norm. If $f$ is $\epsilon_T$-concentrated on $T$ and $\widehat{f}$ is $\epsilon_W$-concentrated on $W,$ then $$|T| \, |W| \geq (1 - (\epsilon_T + \epsilon_W))^2.$$
\end{theorem}

This uncertainty principle has important applications in signal theory and was later generalized to locally compact abelian groups by Smith \cite{Smith-uncertainty_on_lct}. Since its introduction, the Donoho-Stark principle has inspired a large body of research extending its framework to various settings, such as the Hankel transform \cite{Tuan-Dohono_Star_for_Hankel} and the Laguerre-Bessel transform \cite{Negzaoui-Laguerre_transform_benedicks}, where the notion of concentration is adapted to different weighted spaces and measures. For a more detailed discussion of the inequality, one can refer to the references cited in these articles.

\subsection{Uncertainty principle of Price}
In \(\mathbb{R}^n\), Heisenberg’s inequality states that if \(f\) is highly localized, its Fourier transform cannot be concentrated near a single point. The local uncertainty principle extends this classical result by showing that such concentration is impossible even when \(\widehat{f}\) is allowed to occupy small neighborhoods around two or more widely separated points. These inequalities were first studied by Faris \cite{Faris-inequalities_local_uuncertainty} and were later sharpened and generalized by Price \cite{Price-inequalities_local_uuncertainty, Price-sharp_local_uncertainty} and Price in collaboration with Sitaram \cite{Price-Sitaram-Local_uncertainty}. 

Price’s local uncertainty principle was later generalized to various mathematical frameworks, including certain Lie groups \cite{Sitaram-Sundari-Thangavelu-uncertainty_principles}, stratified nilpotent groups \cite{Ciatti-Ricci-Sundari-uncertainity_on_stratified} and more recently, the Laguerre hypergroup \cite{Atef-inequalities_on_laguerre_hypergroup}.

\subsection{Beurling's theorem} 

Among the various uncertainty principles is Beurling's theorem on \(\mathbb{R}^n\), which states that there is no nontrivial function satisfying
$$
\int_{\mathbb{R}^n}{\int_{\mathbb{R}^n} \lvert f(y) \rvert \, \lvert \widehat{f}(\xi) \rvert \, e^{\lvert \langle y,\xi \rangle \rvert}}\,dy\,d\xi < \infty.
$$
 
 Several other uncertainty principles, such
 as the theorems of Hardy, Cowling and Price, and Gelfand and Shilov, can be derived from Beurling's theorem \cite{sarkarbeurling}. The version described above is one of the several formulations of Beurling's theorem.
 
The above version was generalized by Bonami et al. \cite{bonami} as follows:
\begin{theorem}
Let \(f \in L^2(\mathbb{R}^n)\) and \(N \geq 0\). Then
$$
\int_{\mathbb{R}^n }\int_{\mathbb{R}^n } \frac{|f(x)| \, |\widehat{f}(y)|}{(1+\|x\|+\|y\|)^N} \, e^{| \langle x,y \rangle|} \, dx \, dy < \infty 
$$
if and only if $f$ can be written as
$$
f(x) = P(x) \, e^{-\pi \langle Ax,x \rangle},
$$
where \(A\) is a real positive definite symmetric matrix and \(P\) is a polynomial of degree $d$ such that $d < (N-n)/2$.
Also, if $N \leq n$, then $f$ is identically $0$.
\end{theorem}
Elloumi et al. \cite{generalized-beurling} further generalized the version of Beurling's theorem established by Bonami et al. to the case where the integrand is of the form \(|f(x)|^p\,|\widehat{f}(y)|^q\), with \(p,q \geq 1\).
 
Over the past two decades, several analogues of Beurling’s theorem have been established in different mathematical frameworks. In 2007, Huang and Liu \cite{huang-liu} proved an analogue of the theorem for the Laguerre hypergroup, which is the fundamental manifold of the radial function space for the Heisenberg group. In 2008, Parui and Sarkar \cite{Beurling-rudra-parui} established another version for two-step nilpotent Lie groups. More recently, in 2022, Thangavelu \cite{thangavelubeurling} proved Beurling's theorem on the Heisenberg group. H-type groups form a natural generalization of the Heisenberg group, and the second and third authors, in collaboration with Pusti and Thangavelu \cite{Dasgupta-Gulia-Beurling_theorem}, established corresponding analogues of this uncertainty principle within that framework.\\

The study of uncertainty principles within different mathematical frameworks has fascinated researchers for decades. In 2023, Thangavelu \cite{Strichartz_Fourier_thangavelu} introduced a formulation of the Fourier transform on the Heisenberg group, referred to as the Strichartz Fourier transform. This transform is scalar-valued and constructed using the joint eigenfunctions of the sublaplacian \(\mathcal{L}\) and the vector field \(T = \frac{\partial}{\partial t},\) both intrinsic to the structure of the Heisenberg group. It is well known that the classical Fourier transform on the Heisenberg group is operator-valued. Although it satisfies the inversion formula and Plancherel identity, its inherent complexity restricts its applicability to problems such as the characterization of the Schwartz space and the formulation of weighted inequalities. The Strichartz Fourier transform overcomes these challenges by providing a scalar-valued analogue that preserves many structural properties akin to the Helgason Fourier transform on symmetric spaces.

Motivated by this development, the present work undertakes the first systematic investigation of uncertainty principles for the Strichartz Fourier transform on the Heisenberg group. We establish analogues of Benedicks’ theorem, the Donoho-Stark principle, Price’s local uncertainty inequality, and a weak Beurling-type theorem within this framework. A significant contribution is the derivation of a Benedicks-Nazarov uncertainty inequality that quantifies the impossibility of simultaneous localization in both domains. In doing so, we overcome the technical limitations of the traditional operator-valued approach and demonstrate how the scalar-valued Strichartz framework facilitates direct extensions of classical Euclidean results to the noncommutative setting of the Heisenberg group. These results not only highlight the analytical strength of the Strichartz Fourier transform but also lay the foundation for further developments on H-type and reduced H-type groups.

This paper is organized as follows. Section \ref{sec 2} provides the necessary background on the Strichartz Fourier transform and the associated dilation structure on the Heisenberg group. In Section \ref{sec 3}, we establish an analogue of Benedicks’ theorem, from which Nazarov’s uncertainty principle follows as a corollary. Sections \ref{sec 4} and \ref{sec 5} are devoted to the study of the Donoho-Stark uncertainty principle and Price’s local uncertainty inequality, respectively. Finally, Section \ref{sec 6} concludes the paper with a weak form of Beurling’s theorem for the Strichartz Fourier transform.

Throughout this article, the notation $A \lesssim B$ means that there exists a constant $c>0,$ depending only on the parameters, such that $A \leq cB$. Further, we use the notation \(X \asymp Y\) to denote that \(X=cY\) for some constant \(c>0\).

%---------------------------------------------------------------------------

\section{Preliminaries}\label{sec 2}

In this section, we recall essential preliminaries on the Heisenberg group and Laguerre functions, which play a fundamental role in the formulation of the Strichartz Fourier transform. We also review the corresponding inversion formula and Plancherel identity associated with this transform.

Let \(\mathbb{H}^n=\mathbb{C}^n \times \mathbb{R}\) be the \((2n+1)\) dimensional Heisenberg group equipped with the group law 
$$(z,t)(w,s)=\left(z+w,\;t+s+\frac{1}{2}\operatorname{Im}(z\overline{w}) \right),\;\forall \, (z,t),(w,s) \in \mathbb{H}^n.$$
The Heisenberg group is a two-step nilpotent and noncommutative Lie group with a well-developed representation theory. For a detailed exposition of the representation theory of $\mathbb{H}^n$ see \cite{fisher-ruzhansky-book, folland-harmonic_phase_space, rottensteiner-thesis, rottensteiner-ruzhansky-harmonic-oscillators, thangaveluheisenberg} . 

\noindent We denote by \(\mathfrak{h}_n,\) the Heisenberg algebra whose basis consists of the left invariant vector fields 
\begin{equation*}
    X_i=\frac{\partial}{\partial x_i}+\frac{1}{2}y_i\frac{\partial}{\partial t},\; Y_i=\frac{\partial}{\partial y_i}-\frac{1}{2}x_i\frac{\partial}{\partial t},\; T= \frac{\partial}{\partial t},\; i=1,2,\dots , n.
\end{equation*}
The sublaplacian on \(\mathbb{H}^n\) is defined as 
$$
\mathcal{L}:=-\sum_{i=1}^n\left( X_i^2+Y_i^2 \right)
$$
and can be explicitly written as 
$$
\mathcal{L} = -\Delta_{\mathbb{C}^n}
- \frac{1}{4} |z|^2 \frac{\partial^2}{\partial t^2}
+ N \frac{\partial}{\partial t},
$$
where $\Delta_{\mathbb{C}^n}$ denotes the Laplacian on $\mathbb{C}^n$ and
$$
N = \sum_{j=1}^n \left( x_j \frac{\partial}{\partial y_j}
- y_j \frac{\partial}{\partial x_j} \right).
$$ 
It is well known that the sublaplacian $\mathcal{L}$ commutes with the vector field $T$, providing a natural framework for studying their joint spectral theory. The Laguerre functions play a fundamental role in describing the joint eigenfunctions of $\mathcal{L}$ on the Heisenberg group $\mathbb{H}^n$ and the left-invariant vector field $\frac{\partial}{\partial t}$.

For \(\delta >-1\) and \(k \in \mathbb{N} \cup \{0\}\), let \(L^{\delta}_k\) denote the Laguerre polynomial of type \(\delta\) and degree \(k\), given by
$$
L^{\delta}_k(t):=\frac{t^{-\delta}e^t}{k!}\left ( \frac{d}{dt} \right )^k \left( e^{-t} \, t^{k+\delta}\right),\;\;\; t\geq 0.
$$
Let $\mathbb{R}^\ast = \mathbb{R}\setminus\{0\}.$ For \(\lambda \in \mathbb{R}^\ast\), the Laguerre function on \(\mathbb{C}^n\) is denoted by \(\varphi^{n-1}_{k,\lambda}\) and defined as
$$
\varphi^{n-1}_{k,\lambda}(z):=L^{n-1}_k\left( \frac{1}{2}|\lambda||z|^2\right) e^{-\frac{1}{4}|\lambda||z|^2}.
$$
Further details on the Laguerre functions can be found in \cite{beals-gaveau-greiner, beals1988calculus, thangaveluheisenberg, wong_weyl-transforms}.

The study of eigenfunctions of the sublaplacian has been discussed in detail in \cite{thangavelu-uncertainty-book}. 
The functions \(e^{n-1}_{k,\lambda}(z,t):=e^{i \lambda t}\varphi^{n-1}_{k,\lambda}(z)\) satisfy the following eigenvalue equations
$$
\mathcal{L}(e^{n-1}_{k,\lambda})=(2k+n)\lvert \lambda \rvert \,e^{n-1}_{k,\lambda} \;\;\;\; \text{and} \;\;\;-i\frac{\partial}{\partial t}(e^{n-1}_{k,\lambda})=\la \, e^{n-1}_{k,\lambda} .
$$
For \(\lambda \in \mathbb{R}^\ast\), the operators \(L_{\lambda}\) defined by the relation \(\mathcal{L}(e^{i \lambda t}f)(z)=e^{i\lambda t}(L_{\lambda}f)(z)\), are called special Hermite operators. The Laguerre functions \(\varphi^{n-1}_{k,\lambda}\) are eigenfunctions of \(L_\lambda\) with eigenvalue \((2k+n)|\lambda|\) and every function \(f\in L^2(\mathbb{R}^n)\) admits a series representation of the form
$$
f(z)=(2\pi)^{-n}|\lambda|\sum_{k=0}^{\infty}\left(f \ast_{\lambda}\varphi^{n-1}_{k,\lambda}(z)\right).
$$ 
Let \( \alpha \in \mathbb{N}^n \) and denote the Hermite function on \( \mathbb{R}^n \) by \( \Phi_{\alpha} \). It is well known that the collection \( \{\Phi_{\alpha}\}_{\alpha \in \mathbb{N}^n} \) forms an orthonormal basis for \(L^2(\mathbb{R}^n)\), and that each \( \Phi_{\alpha} \) is an eigenfunction of the Hermite operator
\[
H = -\Delta + |x|^2 .
\]
For \( \lambda \in \mathbb{R}^{\ast} \), define the scaled Hermite functions
\[
\Phi_{\alpha}^{\lambda}(x) = |\lambda|^{n/4} \, \Phi_{\alpha}\big(\sqrt{|\lambda|}\,x\big), 
\qquad x \in \mathbb{R}^n .
\]
These functions arise naturally in the context of the Schr\"odinger representation \( \pi_{\lambda} \) of the Heisenberg group \( \mathbb{H}^n \). For \( \alpha, \beta \in \mathbb{N}^n \) with \( |\alpha| = k \), define
\[
\Phi_{\alpha,\beta}^{\lambda}(z) 
= (2\pi)^{-n/2} \, 
\langle \pi_{\lambda}(z,0)\Phi_{\alpha}^{\lambda}, \Phi_{\beta}^{\lambda} \rangle_{L^2(\mathbb{R}^n)}, 
\qquad z \in \mathbb{C}^n \simeq \mathbb{R}^{2n}.
\]
Then the family \( \{\Phi_{\alpha,\beta}^{\lambda} : |\alpha| = k \} \) forms an orthogonal basis for the eigenspace of the Hermite operator \( L_{\lambda} \) associated with the eigenvalue \( (2k+n)|\lambda| \).

As mentioned in the introduction, the classical Fourier transform on 
\(\mathbb{H}^{n}\) is operator--valued, which presents technical difficulties in many standard problems of harmonic analysis. To overcome this, the Strichartz Fourier transform was introduced, formulated in terms of Laguerre and Bessel functions.

The \emph{Heisenberg fan}, denoted by \(\Omega\), is the union of rays
\[
R_{k}
:= \big\{(\lambda,(2k+n)|\lambda|)\,:\, \lambda \in \mathbb{R}^{\ast}\big\},
\qquad k = 0,1,2,\ldots,
\]
together with the limiting ray
\[
R_{\infty}
:= \big\{(0,\tau): \tau \ge 0 \big\}.
\]

According to \cite[Chapter 3]{thangaveluheisenberg}, the \(U(n)\)-spherical functions on \(\mathbb{H}^{n}\) fall into two families:
\[
\frac{k!\,(n-1)!}{(k+n-1)!}
\, e^{\,i\lambda t}\,
\varphi^{\,n-1}_{k,\lambda}(z),
\qquad \text{and} \qquad
2^{\,n-1}(n-1)!
\,\frac{J_{n-1}(\sqrt{\tau}\,|z|)}{(\sqrt{\tau}\,|z|)^{n-1}},
\]
where \(J_{n-1}\) denotes the Bessel function of order \(n-1\), and 
\(\varphi^{\,n-1}_{k,\lambda}\) denotes the Laguerre function of type \(n-1\).

Following Thangavelu \cite{Strichartz_Fourier_thangavelu}, for 
\( f \in L^{1}(\mathbb{H}^{n}) \cap L^{2}(\mathbb{H}^{n}) \), 
the \emph{Strichartz Fourier transform} of \(f\) is defined on 
\(\Omega \times \mathbb{C}^{n}\) by
\begin{equation}\label{sft def1}
    \widehat{f}(a,w)
    =
    \int_{\mathbb{H}^{n}}
        f(z,t)\,
        e_{a}\!\big((z,t)^{-1}(w,0)\big)
    \, dz\,dt,
\end{equation}
where \( a = (\lambda,(2k+n)|\lambda|) \in R_{k} \). 
Throughout, we write
\[
 e_{a}(z,t) := e^{\,n-1}_{k,\lambda}(z,t),
\qquad a = (\lambda,(2k+n)|\lambda|) \in R_{k}.
\]

\noindent For an element \((0,\tau) \in R_{\infty}\), the Strichartz Fourier transform of \(f\) is given by  
\begin{equation*}
    \widehat{f}(0,\tau,w)
    = 2^{n-1}(n-1)! 
      \int_{\mathbb{H}^{n}} 
        f(z,t)\,
        \frac{J_{n-1}\!\left(\sqrt{\tau}\,|w-z|\right)}
             {(\sqrt{\tau}\,|w-z|)^{\,n-1}}\,dz\,dt.
\end{equation*}
This transform provides a convenient framework for the joint spectral 
analysis of the sublaplacian and the central derivative on the 
Heisenberg group.

Since \(\Omega \subseteq \mathbb{R}^{2}\), it inherits the standard Euclidean topology.  
On the set \(\Omega\), define a measure \(\nu\) by  
\begin{equation}\label{nu_def}
    \int_{\Omega} \phi(a)\, d\nu(a)
    := (2\pi)^{-2n-1}
       \int_{-\infty}^{\infty}
          \left(
             \sum_{k=0}^{\infty}
             \phi\big(\lambda,(2k+n)|\lambda|\big)
          \right)
          |\lambda|^{2n}\, d\lambda,
\end{equation}
for a suitable function $\phi$ on $\Omega.$

We now recall the inversion formula and Plancherel identity for the Strichartz Fourier transform.
\begin{theorem}\cite[Theorem 7.2]{Strichartz_Fourier_thangavelu}
    For \(f \in \mathcal{S}(\mathbb{H}^n)\), the inversion formula for the Strichartz Fourier transform is given by 
    \begin{equation}
    f(z,t)=\int_{\Omega \times \mathbb{C}^n}\widehat{f}(a,w) \; e_a \left((-w,0)(z,t)\right) dw \, d\nu(a).
    \end{equation}
Furthermore, for a function \(f\) in \((L^1 \cap L^2)(\mathbb{H}^n)\), the following Plancherel identity holds:

\begin{equation}
    \int_{\mathbb{H}^n}|f(z,t)|^2 \, dz\,dt=\int_{\Omega \times \mathbb{C}^n}|\widehat{f}(a,w)|^2 \, dw \, d\nu(a).
\end{equation}
\end{theorem}

The \emph{normalized Strichartz Fourier transform} of \(f\), denoted 
\(\mathcal{F}(f)=\tilde{f}\), is defined by  
\begin{equation}\label{norm_sft}
    (\mathcal{F}(f))(a,w)
    = \tilde{f}(a,w)
    := c_{n,k}\,\widehat{f}(a,w)
    = \frac{k!\,(n-1)!}{(k+n-1)!}\,\widehat{f}(a,w),
\end{equation}
for \(a=(\lambda,(2k+n)|\lambda|)\in R_{k}\).  
For the degenerate component \((0,\tau)\in R_{\infty}\), we set  
\begin{equation}\label{norm_sft_infty}
    (\mathcal{F}(f))(0,\tau,w)
    = \tilde{f}(0,\tau,w)
    := \widehat{f}(0,\tau,w).
\end{equation}

\noindent Using the estimates for Laguerre and Bessel functions (see \cite{szego})
\[
\frac{k!\,(n-1)!}{(k+n-1)!}\,|\varphi^{\,n-1}_{k,\lambda}(z)| \le 1,
\qquad 
|J_{n-1}(t)| \lesssim t^{n-1},
\]
Thangavelu \cite{Strichartz_Fourier_thangavelu} proved that for any 
\(f \in L^{1}(\mathbb{H}^{n})\), the normalized Strichartz Fourier transform satisfies 
\begin{equation}\label{L1 to L inf bdd}
  \sup_{(a,w)\in \Omega \times \mathbb{C}^{n}}
  \big|(\mathcal{F}(f))(a,w)\big|
  \lesssim \|f\|_{1}.
\end{equation}
Moreover, we have the following inversion formula:
\begin{equation}\label{inversionformula}
   f(z,t)=\int_{\Omega \times \mathbb{C}^n}c_{n,k}\,(\mathcal{F}({f}))(a,w) \; e_a \left((-w,0)(z,t)\right) dw \, d\nu_2(a), 
\end{equation}
and the Plancherel identity takes the form
\begin{equation}\label{plancherel}
    \int_{\mathbb{H}^{n}} |f(z,t)|^{2}\,dz\,dt
    = \int_{\Omega \times \mathbb{C}^{n}}
      \big|(\mathcal{F}(f))(a,w)\big|^{2}\,dw\, d\nu_{2}(a),
\end{equation}
where \(d\nu_{2}(a)\) is a measure on \(\Omega\) given by  
\[
\int_{\Omega} \phi(a)\, d\nu_{2}(a)
:= (2\pi)^{-2n-1}
   \int_{-\infty}^{\infty}
      \left(
           \sum_{k=0}^{\infty}
           \left(\frac{(k+n-1)!}{k!\,(n-1)!}\right)^{2}
           \phi(\lambda,(2k+n)|\lambda|)
      \right)
      |\lambda|^{2n}\,d\lambda.
\]

\noindent In what follows, we use the notation \(\mathcal{F}(f)\) to denote the normalized 
Strichartz Fourier transform, as this convention simplifies the formulation of the 
projection operators introduced in Section \ref{sec 3}.

Recall that the classical Fourier transform on \(\mathbb{R}^{n}\) is a unitary operator from 
\(L^{2}(\mathbb{R}^{n})\) onto itself. 
In contrast, the Strichartz Fourier transform does not define a unitary operator from 
\(L^{2}(\mathbb{H}^{n})\) onto 
\(L^{2}(\Omega \times \mathbb{C}^{n},\, d\nu\,dw)\). 
However, if we restrict our attention to the subspace 
\(L_{0}^{2}(\Omega \times \mathbb{C}^{n},\, d\nu\,dw)\) of \(L^{2}(\Omega \times \mathbb{C}^{n},\, d\nu\,dw)\),
consisting of functions satisfying
\begin{equation}\label{subspace_condition}
    (2\pi)^{-n}|\lambda|^{n}
    \left(
        \varphi^{\,n-1}_{k,\lambda} 
        \ast_{\lambda} 
        \big(\mathcal{F}(f)\big)(a,\cdot)
    \right)(z)
    = \big(\mathcal{F}(f)\big)(a,z),
    \qquad a \in R_{k},
\end{equation}
then the following result holds.

\begin{theorem}\cite[Theorem 7.3]{Strichartz_Fourier_thangavelu}
    The Strichartz Fourier transform is a unitary operator from \(L^2(\mathbb{H}^n)\) onto \(L_0^2(\Omega \times \mathbb{C}^n, d\nu \, dw)\).
\end{theorem}

The Heisenberg group \(\mathbb{H}^{n}\) has homogeneous dimension \(Q = 2n + 2\) 
and admits a one-parameter family of dilations \(\{\delta_{r}\}_{r>0}\), 
defined by 
\[
    \delta_{r}(z,t) := (rz,\, r^{2}t), \qquad r > 0.
\]
The \emph{Korányi norm} on \(\mathbb{H}^{n}\) is given by  
\[
    |(z,t)| := \big(|z|^{4} + t^{2}\big)^{1/4},
\]
which is homogeneous of degree \(1\) with respect to the dilations \(\delta_{r}\).
Let \(S_{K}\) denote the unit sphere in \(\mathbb{H}^{n}\) corresponding to the Korányi norm. 
For any integrable function \(f\) on \(\mathbb{H}^{n}\), the following polar decomposition holds:
\begin{equation}\label{polar_decomp_koranyi}
    \int_{\mathbb{H}^{n}} f(z,t)\,dz\,dt
    = \int_{0}^{\infty}
        \left(
            \int_{S_{K}} f(\delta_{r}y)\, d\sigma(y)
        \right)
        r^{Q-1}\,dr,
\end{equation}
where \(\sigma\) denotes the Radon measure on \(S_{K}\) (see \cite{Pritam-Lacunary_spherical_function}).  

For a function \(f\) on \(\mathbb{H}^{n}\), the dilation of \(f\) is defined by  
\[
    (\delta_{r}f)(z,t) := f(\delta_{r}(z,t)).
\]

The next lemma records the scaling behavior of the normalized Strichartz Fourier transform under group dilations on \(\mathbb{H}^{n}\).  
This will be used to induce a corresponding family of dilations on the parameter space \(\Omega \times \mathbb{C}^{n}\).

\begin{lemma}\label{lem:dilationonfourierside}
Let \(f \in L^{1}(\mathbb{H}^{n})\) and \(r>0\).  
Then the normalized Strichartz Fourier transform satisfies
\begin{equation*}
    \mathcal{F}(\delta_{r}f)(a,w)
    = r^{-Q}\,
      \F(f)\left((r^{-2} \la,(2k+n)|r^{-2} \la|),rw \right),
\end{equation*}
for \(a=(\lambda,(2k+n)|\lambda|)\in R_{k}\).  
For the degenerate component \((0,\tau)\in R_{\infty}\), we have
\begin{equation*}
    \mathcal{F}(\delta_{r}f)(0,\tau,w)
    = r^{-Q}\,
      \mathcal{F}(f)\!\left(0,\,r^{-2}\tau,\,rw\right).
\end{equation*}
\end{lemma}

 \begin{proof}
        Let \(a=(\lambda,(2k+n)|\lambda|)\in R_{k}\) and \(w \in \C^{n}\).  
By the definition of \(\mathcal{F}\), we have
        \begin{align*}
        &\F\left( \delta_r(f)\right)(a,w)\\
        &=c_{n,k} \int_{\Hg^n}(\delta_rf)(z,t) \; e_a\left((z,t)^{-1}(w,0) \right) dz\,dt\\
        &= c_{n,k}\int_{\Hg^n}f(rz,r^2t) \; \varphi^{n-1}_{k,\la}(w-z) \; e^{i\la \left(-t-\frac{1}{2}\operatorname{Im}(z\overline{w}) \right)} \; dz\,dt.
        \end{align*}
        Making the change of variables \(z \mapsto r^{-1}z,\ t \mapsto r^{-2}t\) yields
\begin{align*}
&\F(\delta_{r}f)(a,w)\\
        &=r^{-2n-2} \; c_{n,k}\int_{\Hg^n}f(z,t) \; \varphi^{n-1}_{k,\la}(w-r^{-1}z) \; e^{i\la\left(-r^{-2}t-\frac{1}{2}\operatorname{Im}(r^{-1}z\overline{w}) \right)} \; dz\,dt\\
        &=r^{-Q} \; c_{n,k}\int_{\Hg^n}f(z,t) \; L^{n-1}_k\left(\frac{1}{2}|\la||w-r^{-1}z|^2 \right)e^{\left(-\frac{1}{4}|\la||w-r^{-1}z|^2\right)} \; e^{i\la\left(-r^{-2}t-\frac{1}{2}\operatorname{Im}(r^{-1}z\overline{w}) \right)} \; dz\,dt\\
        &=r^{-Q} \; c_{n,k} \int_{\Hg^n}f(z,t) \; L^{n-1}_k\left(\frac{1}{2}|r^{-2}\la||rw-z|^2 \right) e^{\left(-\frac{1}{4}|r^{-2}\la||rw-z|^2\right)} \; e^{ir^{-2}\la\left(-t-\frac{1}{2}\operatorname{Im}(z\overline{rw}) \right)} \; dz\,dt\\
        &=r^{-Q} \; c_{n,k}\int_{\Hg^n}f(z,t) \; \varphi^{n-1}_{k,r^{-2}\la}(rw-z) \; e^{ir^{-2}\la\left(-t-\frac{1}{2}\operatorname{Im}(z\overline{rw}) \right)} \; dz\,dt\\
        &=r^{-Q} \; c_{n,k} \int_{\Hg^n}f(z,t) \; e_{(r^{-2}\la,(2k+n)|r^{-2}\la|)}\left((z,t)^{-1}(rw,0) \right)dz\,dt\\
        &=r^{-Q} \; \F(f)\left((r^{-2} \la,(2k+n)|r^{-2} \la|),rw \right)
        \end{align*} which proves the result.
    \end{proof}
Motivated by Lemma \ref{lem:dilationonfourierside}, for \(r>0,\) we define the dilation \(\mathcal{D}_r\) on \(\Omega \times \C^n\) by
$$
\mathcal{D}_r(a,w):=\begin{cases}
   \left((r^{2} \la,(2k+n)|r^{2} \la|),r^{-1}w \right), & \text{ if } a=(\la,(2k+n)|\la|)\in R_k\\
   (0,r^2\tau,r^{-1}w),& \text{ if } a = (0,\tau)\in R_{\infty}.
\end{cases}
$$
Thus, we have the following relation:
\begin{equation}\label{dilation-relation}
    \F\left( \delta_r(f)\right) = r^{-Q} \left(\mathcal{D}_{r^{-1}}\F(f)\right),
\end{equation}
where $$\left(\mathcal{D}_{r^{-1}}\F(f)\right)(a,w) := \left( \F(f)\right)(\mathcal{D}_{r^{-1}}(a,w)).$$

In the setting of the Heisenberg group \(\mathbb{H}^{n}\), the Haar measure coincides with the Lebesgue measure on \(\R^{2n+1}\), and hence inherits the familiar regularity properties of Lebesgue measure. More generally, when working on topological groups or parameter spaces such as \(\Omega \times \C^{n}\), it is essential to guarantee that the underlying measure behaves well with respect to the topology. In particular, one would like measurable sets to be approximated by open or compact sets, and integrable functions to be approximated by continuous functions with compact support. This type of regularity is encoded in the notion of a \emph{Radon measure}. The next two classical results ensure that, under mild topological assumptions, Borel measures are automatically Radon and that continuous, compactly supported functions form a dense subspace of \(L^{p}\). These facts will be used repeatedly in our analysis.
The following two results provide sufficient conditions for a measure to be Radon and formalize the density properties of measurable sets and functions in measure spaces equipped with a Radon measure.
\begin{theorem}\cite[Theorem 7.8]{folland-real_analysis_book}\label{measure is radon}
    Let \(X\) be a second-countable, locally compact, Hausdorff space. Then every Borel measure that is finite on compact sets is Radon.
\end{theorem}
\begin{theorem}\cite[Proposition 7.9]{folland-real_analysis_book}\label{density of Cc}
    Let \((X,\mu)\) be a measure space where \(\mu\) is Radon. Then \(\mathcal{C}_c(X)\) is dense in \(L^p(X)\) for all \(1\leq p < \infty\).
\end{theorem}

%-----------------------------------------------------------------------------

\section{Benedicks' Theorem}\label{sec 3}
In this section, we establish an analogue of the Benedicks theorem for the Strichartz Fourier transform, employing projection operators and the notion of the norm of an integral operator.

\medskip
\noindent
Let \(V \subseteq \mathbb{H}^{n}\) and \(W \subseteq \Omega\) be measurable sets, where
\[
W \subseteq \big\{(\lambda,(2k+n)|\lambda|) : |\lambda| > M,\; k \in \mathbb{N} \cup \{0\}\big\}
\]
for some \(M > 0\).  
We consider the following orthogonal projections on \(L^{2}(\mathbb{H}^{n})\):
\[
P_{V}(f) := \chi_{V} f, \qquad 
P_{W}(f) := \mathcal{F}^{-1}\!\left( \chi_{W \times \C^{n}}\, \mathcal{F}(f)\right),
\qquad f \in L^{2}(\mathbb{H}^{n}).
\]

\noindent The operator \(P_{V}\) is clearly well defined.  
For \(P_{W}\) to be well defined, however, we must verify that
\[
\chi_{W \times \C^{n}}\,\mathcal{F}(f) \in 
L^{2}_{0}\!\left(\Omega \times \C^{n},\, d\nu_{2}\,dw\right),
\]
since \(\mathcal{F}\) maps \(L^{2}(\mathbb{H}^{n})\) surjectively onto 
\(L^{2}_{0}\!\left(\Omega \times \C^{n}, d\nu_{2}\,dw\right)\).
In other words, for each \(a = (\lambda,(2k+n)|\lambda|) \in R_{k}\), we require that
\[
(2\pi)^{-n} |\lambda|^{n}
\big(\varphi_{k,\lambda}^{\,n-1} *_{\lambda}
(\chi_{W \times \C^{n}} \mathcal{F}(f))(a,\cdot)\big)(z)
=
(\chi_{W \times \C^{n}} \mathcal{F}(f))(a,z),
\]
where \(*_{\lambda}\) denotes the \(\lambda\)-twisted convolution.

Observe that the above identity holds trivially whenever \(a \notin W\).  
For \(a \in W\), we compute
\begin{align*}
    &(2\pi)^{-n}|\lambda|^{n}
      \Big(\varphi^{\,n-1}_{k,\lambda} *_{\lambda}
      \big(\chi_{W \times \C^{n}}\,\mathcal{F}(f)\big)(a,\cdot)\Big)(z) \\
    &= (2\pi)^{-n}|\lambda|^{n}
       \int_{\C^{n}} 
       \varphi^{\,n-1}_{k,\lambda}(z-w)\,
       \big(\mathcal{F}(f)\big)(a,w)\,
       e^{\,i\frac{\lambda}{2}\Im(z\overline{w})}\,dw \\
    &= (2\pi)^{-n}|\lambda|^{n}
       \big( \varphi^{\,n-1}_{k,\lambda} *_{\lambda} \big(\mathcal{F}(f)\big)(a,\cdot) \big)(z) \\
    &= \big(\mathcal{F}(f)\big)(a,z) \\
    &= \big(\chi_{W \times \C^{n}}\mathcal{F}(f)\big)(a,z),
\end{align*}
where in the second step we used the assumption \(a \in W\), and the third equality follows from the characterization of \(\mathcal{F}(f)\) as it lies in \(L^{2}_{0}(\Omega \times \C^{n},d\nu_{2}\,dw)\).

\noindent Moreover, we clearly have
\[
\|\chi_{W \times \C^{n}} \mathcal{F}(f)\|_{2}
\le \|\mathcal{F}(f)\|_{2}.
\]
Hence \(\chi_{W \times \C^{n}}\mathcal{F}(f)\in 
L^{2}_{0}(\Omega \times \C^{n},d\nu_{2}\,dw)\), and therefore \(P_{W}\) is well defined.

Next we state and prove the following lemma which provides a significant relation between the Hilbert-Schmidt norm of \(P_W \circ P_V\) and the measure of the associated sets \(V\) and \(W\).

\begin{lemma}\label{projection map lemma}
    Let \(V \subseteq \mathbb{H}^n\) and \(W \subseteq \Omega\) such that \(W \subseteq  \left \{(\la,(2k+n)|\la|): |\la|>M,\; k \in \mathbb{N} \cup \{0\} \right \}\), for some \(M>0\). If \(V\) and \(W\) are sets of finite measure, then 
    $$
    \lVert P_W \circ P_V \rVert_{HS}^2 \leq (2\pi)^n M^{-n}|V||W|.
    $$
\end{lemma}

\begin{proof}

For \(f \in L^{2}(\mathbb{H}^{n})\), it follows from the inversion formula \eqref{inversionformula} that
\begin{align}
    &\big((P_{W} \circ P_{V})f\big)(z,t) \nonumber\\
    &= \int_{\Omega \times \C^{n}}
        c_{n,k}\,
        \big(\chi_{W \times \C^{n}}\,\mathcal{F}(P_{V}f)\big)(a,w)\;
        e_{a}\!\left((-w,0)(z,t)\right)\,
        d\nu_{2}(a)\,dw \nonumber \\
    &= \int_{\Omega \times \C^{n}}
        c_{n,k}\,\chi_{W \times \C^{n}}(a,w)
        \left(
            \int_{\mathbb{H}^{n}}
                c_{n,k}\,(P_{V}f)(z',t')\,
                e_{a}\!\left((z',t')^{-1}(w,0)\right)
            dz'\,dt'
        \right) \nonumber \\
        &\hspace{10cm}e_{a}\!\left((-w,0)(z,t)\right)\,
        d\nu_{2}(a)\,dw.
        \label{expression of PwPv}
\end{align}

For \((z',t') \in \mathbb{H}^{n}\), define a function \(h_{(z',t')} : \Omega \times \C^{n} \to \C\) by
\[
    h_{(z',t')}(a,w):=
    \begin{cases}
        c_{n,k}\,\chi_{W \times \C^{n}}(a,w)\,
        e_{a}\!\left((z',t')^{-1}(w,0)\right), & a \in R_{k}, \\[4pt]
        0, & a \in R_{\infty}.
    \end{cases}
\]
We claim that \(h_{(z',t')} \in L^{2}_{0}(\Omega \times \C^{n},d\nu_{2}\,dw)\). Indeed,
\begin{align}
    \|h_{(z',t')}\|_{2}^{2}
    &= (2\pi)^{-2n-1}
       \int_{\C^{n}} \int_{-\infty}^{\infty}
       \sum_{k=0}^{\infty}
       (c_{n,k})^{-2}
       \left| h_{(z',t')}\left((\la,(2k+n)| \la|),w \right) \right|^{2}
       |\lambda|^{2n}
       d\lambda\,dw \nonumber\\
    &= (2\pi)^{-2n-1}
       \int_{\C^{n}}\int_{-\infty}^{\infty}
       \sum_{k=0}^{\infty}
       (c_{n,k})^{-2}
       \left| c_{n,k} \, \chi_{W \times \C^{n}}\left((\la,(2k+n)| \la|),w \right)
       \varphi^{\,n-1}_{k,\lambda}(w-z') \right|^{2} \nonumber \\
       &\hspace{12.5cm }|\lambda|^{2n}
       d\lambda\,dw \nonumber \\
    &= (2\pi)^{-2n-1}
       \int_{-\infty}^{\infty}
       \sum_{k=0}^{\infty}
       \chi_{W}(\lambda,(2k+n)|\lambda|)
       \left(
           \int_{\C^{n}}
               \left|\varphi^{\,n-1}_{k,\lambda}(w-z')\right|^{2}
           dw
       \right)
       |\lambda|^{2n}
       d\lambda \nonumber\\
    &= (2\pi)^{n}(2\pi)^{-2n-1}
       \int_{-\infty}^{\infty}
       \sum_{k=0}^{\infty}
       (c_{n,k})^{-2} \, \chi_{W}(\lambda,(2k+n)|\lambda|)\,
       c_{n,k}\,|\lambda|^{-n}\,|\lambda|^{2n}d\lambda \nonumber\\
    &\le (2\pi)^{n}M^{-n}
       \left(
           (2\pi)^{-2n-1}
           \int_{-\infty}^{\infty}
           \sum_{k=0}^{\infty}
           (c_{n,k})^{-2} \, \chi_{W}(\lambda,(2k+n)|\lambda|)
           |\lambda|^{2n}d\lambda
       \right) \nonumber\\
    &= (2\pi)^{n}M^{-n}|W| < \infty.
    \label{estimate of norm h}
\end{align}
%Therefore \(h_{(z',t')} \in L^{2}_{0}(\Omega \times \C^{n},d\nu_{2}\,dw)\), as claimed.
 The above calculation uses the following relations:
    $$
    |\la|>M \text{ implies } |\la|^{-n}< M^{-n},
    $$
    $$
    c_{n,k} \leq  1  \; \; \forall  \; k,
    $$
    and 
    \begin{align*}
      \int_{\mathbb{C}^n}|\varphi^{n-1}_{k,\la}(w)|^2dw&= (2\pi)^n |\la|^{-n} \sum_{|\al|=k}1 =(2\pi)^n |\la|^{-n}  \; (c_{n,k})^{-1}.
    \end{align*}
Thus, the function \(h_{(z',t')}\in L^2(\Omega \times \mathbb{C}^n,d\nu_2\,dw)\) for each \((z',t')\in \mathbb{H}^n\).

\noindent Further, for \(a\in W\cap R_k\),
\begin{align*}
 &(2\pi)^{-n}|\la|^n \left( \varphi^{n-1}_{k,\la} \ast_\la \left(h_{(z',t')} \right)(a, \cdot) \right)(z)\\
 &=(2\pi)^{-n}|\la|^n \int_{\C}\varphi^{n-1}_{k,\la}(z-w)h_{(z',t')}(a,w)\,e^{i\frac{\la}{2}\operatorname{Im}(z\overline{w})}dw\\
 &=(2\pi)^{-n}|\la|^nc_{n,k} \int_{\C}\varphi^{n-1}_{k,\la}(z-w)\chi_{W \times \mathbb{C}^n}(a,w)  \; e_a\left((z',t')^{-1}(w,0) \right)e^{i\frac{\la}{2}\operatorname{Im}(z\overline{w})}dw\\
& =(2\pi)^{-n}|\la|^nc_{n,k} \int_{\C}\varphi^{n-1}_{k,\la}(z-w)  \; e_a\left((z',t')^{-1}(w,0) \right)e^{i\frac{\la}{2}\operatorname{Im}(z\overline{w})}dw\\
&=(2\pi)^{-n}|\la|^nc_{n,k} \int_{\C}\varphi^{n-1}_{k,\la}(z-w)  \; \varphi^{n-1}_{k,\la}(w-z')  \; e^{i \la(-t'-\frac{1}{2}\operatorname{Im}(z'\overline{w}))}e^{i\frac{\la}{2}\operatorname{Im}(z\overline{w})}dw\\
&=(2\pi)^{-n}|\la|^nc_{n,k}  \; e^{i\frac{\la}{2}\operatorname{Im}(z\overline{z'})}  \; e^{-i\la t'} \int_{\C}\varphi^{n-1}_{k,\la}(z-z'-w)  \; \varphi^{n-1}_{k,\la}(w)  \; e^{i\frac{\la}{2}\operatorname{Im}\left((z-z')\overline{w} \right)}dw\\
&=(2\pi)^{-n}|\la|^nc_{n,k}  \; e^{i\frac{\la}{2}\operatorname{Im}(z\overline{z'})}  \; e^{-i\la t'} \left(\varphi^{n-1}_{k,\la} \ast_\la \varphi^{n-1}_{k,\la}\right)(z-z')\\
&=c_{n,k}  \; e^{i\frac{\la}{2}\operatorname{Im}(z\overline{z'})}  \; e^{-i\la t'}\varphi^{n-1}_{k,\la}(z-z')\\
&=h_{(z',t')}(a,z).
\end{align*}
This implies that \(h_{(z',t')}\in L^2_0(\Omega \times \mathbb{C}^n,d\nu_2\,dw)\), thus proving the above claim.

Now, consider the function $\mathcal{N}$ on $\mathbb{H}^n \times \mathbb{H}^n$ given by 
\begin{equation*}
    \mathcal{N}(z',t',z,t):=\chi_V(z',t') \left( \F^{-1}\left( h_{(z',t')}\right) \right) (z,t).
\end{equation*}
We claim that $P_W \circ P_V$ is an integral operator with kernel $\mathcal{N},$ that is,
\begin{equation*}
    \left( \left(P_W \circ P_V \right)f \right)(z,t)= \int_{\Hg^n}f(z',t')  \; \mathcal{N}(z',t',z,t)\,dz'\,dt', \hspace{1cm} (f \in L^2(\mathbb{H}^n)).
\end{equation*}
To establish this claim, we explicitly compute the integral on the RHS. Observe that 
\begin{align}
& \int_{\Hg^n}f(z',t')  \; \mathcal{N}(z',t',z,t)\,dz'dt' \nonumber \\
&= \int_{\Hg^n}f(z',t') \; \chi_V(z',t') \; \left( \F^{-1}\left( h_{(z',t')}\right)\right)(z,t) \; dz'dt' \nonumber \\
&= \int_{\Hg^n} \left( P_Vf \right)(z',t') \int_{\Omega \times \C^n}c_{n,k} \; h_{(z',t')}(a,w) \; e_a\left((-w,0)(z,t) \right)\,d\nu_2(a)\,dw\,dz'dt' \nonumber \\
&= \left( \left(P_W \circ P_V \right)f \right) (z,t),\nonumber
\end{align}
where the last equality follows from \eqref{expression of PwPv}.

Therefore, 
\begin{equation*}
    \lVert P_W \circ P_V\rVert^2_{HS}=\lVert \mathcal{N}\rVert_2^2.
\end{equation*}
Now, by Plancherel identity and \eqref{estimate of norm h}, we have
\begin{align*}
    \lVert \mathcal{N}\rVert_2^2&= \int_{\Hg^n \times \Hg^n}\left|\chi_V(z',t') \right|^2 \left|\left(\F^{-1}\left(h_{(z',t')} \right)\right)(z,t) \right|^2dz\,dt\,dz' dt'\\
    &= \int_{\mathbb{H}^n} \left|\chi_V(z',t') \right|^2 \left( \int_{\mathbb{H}^n} \left|\left(\F^{-1}\left(h_{(z',t')} \right)\right)(z,t) \right|^2 dz \, dt \right) dz' dt'\\
    &\leq (2\pi)^n M^{-n}|W| \int_{\mathbb{H}^n} \left|\chi_V(z',t') \right|^2 dz' \, dt' = (2\pi)^n M^{-n} |V| |W|.
\end{align*}
This completes the proof.
\end{proof}
Here we recall the following lemma from \cite{jaming-integral_operators}. Let \(K\) be a convex cone in \(\mathbb{R}^d\) endowed with a suitable Borel measure \(\mu\) which is absolutely continuous with respect to the Lebesgue measure.
\begin{lemma}\cite[Lemma 3.4]{jaming-integral_operators}\label{jaming lemma}
    For a nonzero function \(f \in L^2(K,\mu)\) whose support has finite measure, the family of dilates of \(f\) is linearly independent.
\end{lemma}
Next as a preliminary step, by invoking the above lemma, we establish the following elementary fact about dilations on \(\mathbb{H}^n\).

\begin{prop}\label{linearly independent dilates}
   Let \(f \in L^{2}(\Hg^{n})\) be a nonzero function whose support has finite measure.  
Then the family of dilates \(\{\delta_{r}f : r > 0\}\) is linearly independent. 
\end{prop}

\begin{proof}
    Let \(c_1,c_2,\cdots,c_m\) be scalars such that 
    \begin{equation*}
        \sum_{i=1}^m c_i \, (\delta_{r_i}f)(z,t)=0, \hspace{1cm} r_i > 0.
    \end{equation*}
   Without loss of generality, if we define a function \(g(z):=f(z,0)\), then \(g \in L^2(\C^n)\) is a nonzero function with support of finite measure. Now, by using Lemma \ref{jaming lemma} for the function \(g\), it follows that \(c_i=0\) for all \(i = 1,\cdots,m.\) 
\end{proof}

We recall that the density of \(\mathcal{C}_{c}\)-functions in \(L^{p}\)-spaces \((1 \le p < \infty)\) plays a fundamental role in analysis. Since the Haar measure on \(\Hg^{n}\) coincides with Lebesgue measure, it follows that \(\mathcal{C}_{c}(\Hg^{n})\) is dense in \(L^{p}(\Hg^{n})\).  
To establish an analogous statement on the spectral side, we must examine the topology of \(\Omega\).  

Let \(K \subseteq \Omega\) be compact. As \(\Omega \subseteq \R^{2}\) inherits the Euclidean topology, if \((\lambda,(2k+n)|\lambda|) \in K\), then there exist constants \(B,D>0\) such that \(|\lambda| < B\) and \((2k+n)|\lambda| < D\).  
These bounds imply that the index \(k\) can take only finitely many values, say \(k_{1},k_{2},\dots,k_{m}\). In other words, any compact subset of \(\Omega\) intersects only finitely many rays \(R_{k}\), a fact that will prove useful in establishing approximation results on \(\Omega\).

\noindent Thus, 
    \begin{align}
        |K|&=(2\pi)^{-2n-1}\int_{-\infty}^{\infty}\sum_{k=0}^\infty (c_{n,k})^{-2}\chi_K\left( \lambda,(2k+n)|\lambda|\right)|\lambda|^{2n}d\lambda \nonumber\\
        &=(2\pi)^{-2n-1}\int_{-B}^{B}\sum_{i=1}^m (c_{n,k_i})^{-2}|\lambda|^{2n}d\lambda <\infty. \nonumber 
    \end{align}
    Since every compact subset of \(\Omega\) has finite measure, it follows from Theorem \ref{measure is radon} that the measure on \(\Omega\) is Radon. Consequently, by Theorem \ref{density of Cc}, we conclude that \(\mathcal{C}_{c}(\Omega)\) is dense in \(L^{p}(\Omega)\) for all \(1 \le p < \infty\).\\
We now state the Strichartz Fourier transform version of the Benedicks theorem.

\begin{theorem}\label{Benedicks theorem}
  Let \(V \subseteq \Hg^{n}\) and \(W \subseteq \Omega\) be measurable sets, and assume that
\[
W \subseteq \big\{ (\lambda,(2k+n)|\lambda|) : |\lambda| > M,\; k \in \mathbb{N} \cup \{0\} \big\},
\]
for some \(M>0\). Suppose further that both \(V\) and \(W\) have finite measure.  
If \(f \in L^{2}(\Hg^{n})\) satisfies \(\operatorname{supp}(f) \subseteq V\) and
\[
\operatorname{supp}(\mathcal{F}(f)) \subseteq W \times B,
\qquad \text{for some measurable } B \subseteq \C^{n},
\]
then \(f = 0\).

\end{theorem}

\begin{proof}
    Assume, to the contrary, that there exists a nonzero function \(f_{0} \in L^{2}(\Hg^{n})\) such that 
\(\operatorname{supp}(f_{0}) \subseteq V_{0}\) and 
\(\operatorname{supp}(\mathcal{F}(f_{0})) \subseteq W_{0} \times B_{0}\) for some measurable \(B_{0} \subseteq \C^{n}\), 
where \(V_{0},W_{0}\) have finite measure and
\[
W_{0} \subseteq \big\{ (\lambda,(2k+n)|\lambda|) : |\lambda| > M,\; k \in \mathbb{N} \cup \{0\} \big\},
\]
for some \(M > 0\).  
Choose measurable sets \(V_{1} \subseteq \Hg^{n}\) and \(W_{1} \subseteq \Omega\) such that 
\(V_{0} \subseteq V_{1}\), \(W_{0} \subseteq W_{1}\), and \(|V_{1}|,|W_{1}| < \infty\).

 We now claim that there exists an infinite sequence of distinct numbers \(\{r_{i}\}_{i=0}^{\infty}\) with \(r_{0}=1\) such that
\[
\frac{1}{2} < r_{i} < \frac{3}{2} \quad \text{for all } i,
\]
and
\[
|V| < 2|V_{0}|, \qquad |W| < 2|W_{0}|,
\]
where 
\[
V = \bigcup_{i=0}^{\infty} \delta_{r_{i}^{-1}} V_{0}
\quad \text{and} \quad
W = \bigcup_{i=0}^{\infty} \mathcal{D}_{r_{i}} W_{0}.
\]

    \noindent To construct such a sequence, it suffices to verify the continuity of the following maps:

\begin{enumerate}
    \item \(r \overset{\phi_S}{\longmapsto} |S \cup \delta_rV_0|\), for \(r \in (0,\infty)\) and a subset \(S\) of \(\Hg^n\) of finite measure.
    \item \(r \overset{\psi_T}{\longmapsto} |T\cup\mathcal{D}_r W_0|\), for \(r \in (0,\infty)\) and a finite measure subset \(T\) of \(\Omega\).
    \item \(r \overset{\rho}{\longmapsto} r^{-1} \), for \(r \in (0,\infty)\).
\end{enumerate}

The continuity of \(\rho\) is immediate.  
For the maps \(\phi_{S}\) and \(\psi_{T}\), note that
\begin{eqnarray*}
   & |S \cup \delta_rV_0|=\lVert \chi_{\delta_rV_0}-\chi_{S}\rVert_2^2+\langle \chi_{\delta_rV_0}, \chi_{S}\rangle,\\
   & |T \cup \mathcal{D}_rW_0|=\lVert \chi_{\mathcal{D}_rW_0}-\chi_{T}\rVert_2^2+\langle \chi_{\mathcal{D}_rW_0}, \chi_{T}\rangle.
\end{eqnarray*}
By Theorem \ref{density of Cc}, characteristic functions of finite-measure sets can be approximated in \(L^{2}\) by continuous compactly supported functions.  
Since the dilation maps \(r \mapsto \delta_{r}V_{0}\) and \(r \mapsto \mathcal{D}_{r}W_{0}\) act continuously in measure, it follows that \(\phi_{S}\) and \(\psi_{T}\) are continuous functions of \(r\).

\noindent Now, if we take \(S=V_1\), then by continuity of \(\phi_{V_1}\) at \(1\), there exists an open set \(\mathcal{O}_1\) containing \(1\), such that 
\begin{equation*}
    |V_1 \cup \delta_rV_0|-|V_1 \cup V_0|<\frac{|V_0|}{4},\; \text{ for } r \in \mathcal{O}_1.
\end{equation*}
Similarly, by continuity of \(\psi_{W_1}\) at \(1\), there exists an open set \(\mathcal{O}_2\) containing \(1\), such that 
\begin{equation*}
    |W_1 \cup \mathcal{D}_r W_0|-|W_1 \cup W_0|<\frac{|W_0|}{4},\; \text{ for } r \in \mathcal{O}_2.
\end{equation*}
Further, by continuity of \(\rho\) at \(1\), \(\rho^{-1}(\mathcal{O}_1)\) is an open set containing \(1\).

If we take \(\mathcal{V}_1= \rho^{-1}(\mathcal{O}_1)\cap \mathcal{O}_2 \cap \left(\frac{1}{2},\frac{3}{2} \right)\), then for \(r \in \mathcal{V}_1\), the following holds:
\begin{eqnarray*}
&|V_1 \cup \delta_{{r}^{-1}}V_0|-|V_1|<\frac{|V_0|}{4},\\
&|W_1 \cup \mathcal{D}_{r} W_0|-|W_1|<\frac{|W_0|}{4},\\
& r \in \left(\frac{1}{2},\frac{3}{2} \right) .   
\end{eqnarray*}
Choose an element \(r_1 \neq 1=r_0\) in \( \mathcal{V}_1\) and let \(V_2=\delta_{{r_0}^{-1}}V_0 \cup \delta_{{r_1}^{-1}}V_1 \) and \(W_2=\mathcal{D}_{r_0}W_0 \cup \mathcal{D}_{r_1}W_1 \). It is clear that both \(V_2\) and \(W_2\) are sets of finite measure such that \(V_0 \subseteq V_2\) and \(W_0 \subseteq W_2\). Repeating the above steps, we obtain \(r_2 \in \left(\frac{1}{2},\frac{3}{2}\right)\), distinct from \(r_0\) and \(r_1\) such that the following holds:
\begin{eqnarray*}
&|V_2 \cup \delta_{{r_2}^{-1}}V_0|-|V_2|<\frac{|V_0|}{8},\\
&|W_2 \cup \mathcal{D}_{r_2} W_0|-|W_2|<\frac{|W_0|}{8}.  
\end{eqnarray*}
Following this approach, we construct \(V_n= \bigcup_{i=0}^{n-1}\delta_{{r_i}^{-1}}V_0\) and \(W_n= \bigcup_{i=0}^{n-1}\mathcal{D}_{r_i}W_0\) and choose \(r_n\) such that
\begin{eqnarray*}
    &r_n\neq r_i,\; 0 \leq i \leq n-1,\\
   &|V_n \cup \delta_{{r_n}^{-1}}V_0|-|V_n|<\frac{|V_0|}{2^{n+1}},\\
&|W_n \cup \mathcal{D}_{r_n} W_0|-|W_n|<\frac{|W_0|}{2^{n+1}},\\
& r_n \in \left( \frac{1}{2}, \frac{3}{2}\right).
\end{eqnarray*}

Define 
\[
V = \bigcup_{i=0}^{\infty} \delta_{r_{i}^{-1}}V_{0}
\qquad \text{and} \qquad
W = \bigcup_{i=0}^{\infty} \mathcal{D}_{r_{i}}W_{0}.
\]
A straightforward measure estimate shows that \(|V| < 2|V_{0}|\) and \(|W| < 2|W_{0}|\), thereby completing the argument.

Now, if we define \(f_i:=\delta_{r_i}f_0\), then \(\operatorname{supp}(f_i) \subseteq \delta_{{r_i}^{-1}}V_0 \subseteq V\) for each \(i\). Further, observe that $$\operatorname{supp}\left(\F(f_i) \right)=\operatorname{supp} \left(\F(\delta_{r_i}f_0) \right) \subseteq \mathcal{D}_{r_i} W_0 \times \mathbb{C}^n\subseteq W \times \C^n.$$ Moreover, if \(\left( \la,(2k+n)|\la|\right) \in W\), then \(\left( \la,(2k+n)|\la|\right) \in \mathcal{D}_{r_i}W_0\) for some \(i\). Thus, 
$$
\left( \la,(2k+n)|\la|\right)=\left(r_i^2 \la',(2k+n)|r_i^2 \la'| \right),
$$
where \(\la' \in W_0\). 
This implies that \(|\lambda| = |r_{i}^{2}\lambda'| > \tfrac{1}{4}M\), and therefore the set \(W\) satisfies the hypothesis of Lemma \ref{projection map lemma}. By the definitions of \(P_{V}\) and \(P_{W}\), we have
\[
P_{V}(f_{i}) = f_{i} = P_{W}(f_{i}) \quad \text{for every } i.
\] 
Since the sequence \(\{r_{i}\}_{i=0}^{\infty}\) consists of distinct numbers, Proposition \ref{linearly independent dilates} yields
\[
\dim\big(\operatorname{Im}(P_{V}) \cap \operatorname{Im}(P_{W})\big) = \infty.
\]

\noindent However, by Lemma \ref{projection map lemma}, we have $$
\dim\big(\operatorname{Im}(P_V) \cap \operatorname{Im}(P_W) \big)\leq \lVert P_W \circ P_V \rVert^2_{HS} < \infty,
$$
which is a contradiction. Hence, the result follows.
\end{proof}

As a corollary to Benedicks' theorem, we derive the Nazarov's uncertainty principle.

\begin{corollary}
    Let \(V \subseteq \Hg^{n}\) and \(W \subseteq \Omega\).  
If \(V\) and \(W\) satisfy the hypotheses of Theorem \ref{Benedicks theorem}, then there exists a constant \(C>0\) such that for every \(f \in L^{2}(\Hg^{n})\),
\[
\|f\|_{2}^{2}
\le C\left(
    \|f\|_{L^{2}(V^{c})}^{2}
    + \|\mathcal{F}(f)\|_{L^{2}((W \times B)^{c})}^{2}
\right),
\]
where \(B \subseteq \C^{n}\) is a measurable subset of finite measure.
\end{corollary}

\begin{proof}
If possible, assume that there does not exist a constant \(D\) which satisfies 
    $$
    \lVert f \rVert^2_2\leq D \, \lVert \F(f) \rVert^2_{L^2{\left(({W \times B}\right)^c)}}
    $$
    for all \(f\in L^2(\Hg^n)\) which are supported in \(V\). This implies there exist a sequence \(\{f_n\}\) in \(L^2(\Hg^n)\) with \(\lVert f_n \rVert_2=1\) and supported in \(V\) such that 
    $$
    \lVert \chi_{(W \times B)^c} \, \F (f_n)\rVert_2 \to 0.
    $$
    By Banach-Alaoglu theorem, it follows that the sequence \(\{f_n\}\) converges weakly. Let us call the weak limit to be \(f\). Recall that
    $$
    \left(\F(f_n)\right)(a,w)=c_{n,k}\int_{\Hg^n}f_n(z,t) \; e_a\left((z,t)^{-1}(w,0) \right)\,dz\,dt.
    $$
    Since the function \(\chi_Ve_a \in L^2(\Hg^n) \), it follows from the weak convergence of \(\{f_n\}\) that the sequence \(\left\{\F(f_n)\right\} \) converges to \(\F(f)\).
    Further, observe that $$|\F(f_n)|\leq \int_{\Hg^n}\lvert f_n(z,t) \rvert \,dz\,dt \leq \sqrt{|V|}.$$ Since $B \subseteq \mathbb{C}^n$ is of finite measure, by Lebesgue's dominated convergence theorem, it follows that the sequence \(\left\{\chi_{(W \times B)} \, \F(f_n)\right\}\) converges to \(\chi_{(W \times B)} \, \F(f).\) Moreover, one can easily deduce that \(\operatorname{supp}(f) \subseteq V\) and \(\operatorname{supp}(\F(f)) \subseteq W\times B\). Now, by using Theorem \ref{Benedicks theorem}, we deduce that \(f=0\), which is a contradiction.
\end{proof}

%-----------------------------------------------------------------------------

\section{Donoho-Stark Principle}\label{sec 4}

The Benedicks theorem for the Strichartz Fourier transform asserts that if a function and its Strichartz Fourier transform are supported on sets satisfying the hypotheses of Theorem \ref{Benedicks theorem}, then the function must be identically zero.  
A closely related manifestation of the uncertainty principle is given by the Donoho-Stark theorem \cite{Donoho-stark_paper}, which establishes a similar impossibility of simultaneous localization in mutually dual domains.

We next establish a key lemma, which will be instrumental in proving the Donoho-Stark type estimate on \(\Hg^{n}\).  
Here \(V^{c}\) and \(W^{c}\) denote the complements of \(V\) and \(W\), respectively.

\begin{lemma}
  Let \(V \subseteq \mathbb{H}^n\) and \(W \subseteq \Omega\) such that \(W \subseteq  \left \{(\la,(2k+n)|\la|): |\la|>M,\; k \in \mathbb{N} \cup \{0\} \right \}\), for some \(M>0\). If \(\lVert P_W \circ P_V \rVert<1\), then for all \(f \in L^2(\Hg^n)\),
  \begin{equation}\label{lemma for donoho}
      \lVert f\rVert^2_2\leq \left(1-\lVert P_W \circ P_V \rVert \right)^{-2}\left( \lVert P_{V^c}f \rVert^2_2+\lVert P_{W^c}f \rVert^2_2 \right).
  \end{equation}
  
 \begin{proof}
     Let \(I\) be the identity map on \(L^2(\Hg^n)\). Observe that 
     $$I=P_W+P_{W^c}=P_W(P_V+P_{V^c})+P_{W^c}.$$
     Thus, for \(f \in L^2(\Hg^n)\), we have
     \begin{align*}
         \lVert f-\left(P_W\circ P_V\right)(f)\rVert_2^2&=\lVert \left(P_W \circ P_{V^c}\right)(f)+P_{W^c}f\rVert^2_2\\
         &=\lVert \left(P_W \circ P_{V^c}\right)(f)\rVert_2^2+\lVert P_{W^c}f\rVert^2_2\\
         &\leq \lVert P_{V^c}f\rVert_2^2+\lVert P_{W^c}f\rVert^2_2.
     \end{align*}
     This implies that
     \begin{equation}\label{donoho_lemma_eqn1}
         \lVert f- \left(P_W\circ P_V\right)(f)\rVert_2\leq \left( \lVert P_{V^c}f\rVert_2^2+\lVert P_{W^c}f\rVert^2_2\right)^{1/2}.
     \end{equation}
     Further, by the triangle inequality, we have
     \begin{align*}
         \lVert f- \left(P_W\circ P_V\right)(f)\rVert_2 &\geq \lVert f \rVert_2-\lVert \left(P_W\circ P_V\right)(f) \rVert_2\\
         &\geq \lVert f \rVert_2-\lVert P_W\circ P_V \rVert\lVert f\rVert_2.
     \end{align*}
     Combining this with \eqref{donoho_lemma_eqn1} and using the fact that \(\lVert P_W \circ P_V\rVert<1\), we obtain \eqref{lemma for donoho}.
 \end{proof}
\end{lemma}

 The Donoho-Stark uncertainty principle provides a quantitative lower bound on the product of the measures of the sets on which a function is ``time-limited" and ``band-limited". 
 
\begin{defi}
Let \(V \subseteq \Hg^{n}\) and \(W \subseteq \Omega\) be measurable sets.  
For \(f \in L^{2}(\Hg^{n})\) and parameters \(\epsilon_{V}, \epsilon_{W} > 0\), we say:
\begin{itemize}
    \item \(f\) is \(\epsilon_{V}\)-\emph{time-limited} on \(V\) if 
    \[
    \|P_{V^{c}} f\|_{2} \le \epsilon_{V} \|f\|_{2};
    \]
    \item \(f\) is \(\epsilon_{W}\)-\emph{band-limited} on \(W\) if
    \[
    \|P_{W^{c}} f\|_{2} \le \epsilon_{W} \|f\|_{2}.
    \]
\end{itemize}
\end{defi}

We are now prepared to state the Donoho-Stark uncertainty principle for the Strichartz Fourier transform.

\begin{theorem}
Let \(V \subseteq \Hg^{n}\) and \(W \subseteq \Omega\) be measurable sets, and assume that
\[
W \subseteq \big\{ (\lambda,(2k+n)|\lambda|) : |\lambda| > M,\; k \in \mathbb{N} \cup \{0\} \big\},
\]
for some \(M>0\). Let \(\epsilon_{V}, \epsilon_{W} > 0\) satisfy \(\epsilon_{V}^{2} + \epsilon_{W}^{2} < 1\).  
If \(f \in L^{2}(\Hg^{n})\) is a nonzero function that is both \(\epsilon_{V}\)-time-limited on \(V\) and \(\epsilon_{W}\)-band-limited on \(W\), then
\[
|V|\,|W|
\;\ge\;
(2\pi)^{-n} M^{\,n}
\left(1 - \sqrt{\epsilon_{V}^{2} + \epsilon_{W}^{2}}\right)^{2}.
\]
\end{theorem}

\begin{proof}
If \(|V|\,|W| \ge (2\pi)^{-n} M^{\,n}\), there is nothing to prove.  
Thus assume \(|V|\,|W| < (2\pi)^{-n} M^{\,n}\).  
By Lemma \ref{projection map lemma},
\[
\|P_{W} \circ P_{V}\|
    \;\le\;
    \sqrt{(2\pi)^{n} M^{-n} |V|\,|W|}
    \;<\; 1.
\]
Hence \(V\) and \(W\) satisfy the assumptions of Lemma \ref{lemma for donoho}, and therefore
\[
\big(1 - \|P_{W} \circ P_{V}\|\big)^{-2}
    \;\ge\;
    \frac{\|f\|_{2}^{2}}{\|P_{V^{c}}f\|_{2}^{2} + \|P_{W^{c}}f\|_{2}^{2}}
    \;\ge\;
    \frac{1}{\epsilon_{V}^{2} + \epsilon_{W}^{2}}.
\]
Rearranging yields
\[
1 - \sqrt{\epsilon_{V}^{2} + \epsilon_{W}^{2}}
    \;\le\; \|P_{W} \circ P_{V}\|
    \;\le\;
    \sqrt{(2\pi)^{n} M^{-n} |V|\,|W|},
\]
from which the desired inequality follows.
\end{proof}

%------------------------------------------------------------------------

\section{Uncertainty Principle of Price}\label{sec 5}

In 1986, Price proved the following uncertainty principle for the Fourier transform on \(\R^n\) \cite[Theorem 1.1]{Price-sharp_local_uncertainty}: If \(E\) is a subset of \(\mathbb{R}^n\) and \(\alpha>n/2\), then the following inequality
\begin{equation*}
   \int_{E}| \widehat{f}(\xi)|^2 d\xi \lesssim |E| \; \lVert f \rVert_2^{2-n/\alpha} \lVert \, |x|^{\alpha}f \rVert_2^{n/\alpha}
\end{equation*}
holds for all \(f \in L^2(\mathbb{R}^n)\), where \(\widehat{f}(\xi)\) denotes the Fourier transform of \(f\) in \(\mathbb{R}^n\).

Furthermore, in 1988, Price and Sitaram \cite{Price-Sitaram-Local_uncertainty} established the following inequality
$$
\int_{E}| \widehat{f}(\xi)|^2 \, d\xi \lesssim |E|^{2\alpha/n} \lVert \, |x|^{\alpha}f \rVert_2^{2}
$$
for the case \(0<\alpha<n/2\).

In this section, we establish the local uncertainty principle of Price for the Strichartz Fourier transform. We begin with the following proposition.

\begin{prop}\label{prop sum inequality}
Let \(1 \leq p,q \leq \infty\) and $\alpha, \lambda, \mu > 0.$ Suppose that \(p < q\) and \( \alpha > Q \beta\) where \( \beta = 1/p - 1/q.\) Then for all \(f \in L^q(\mathbb{H}^n)\), the following inequality holds:
$$\|f\|_p^q \lesssim \; \lambda^{-1} \; (\lambda/\mu)^{Q \beta/\alpha} \; \left( \lambda \, \|f\|_q^q + \mu \, \| \, |(z,t)|^\alpha f\|_q^q \right).$$
\end{prop}

\begin{proof}
    Let \( s=q/p.\) By the given condition on \(p\) and \(q\), it follows that \(s > 1.\) Now, using H\"older's inequality, we have
    \begin{align*}
    \lVert f \rVert_p^q&= \left(\int_{\mathbb{H}^n}|f|^p \, dz\,dt \right)^s\\
    &\leq\left( \left\{ \int_{\mathbb{H}^n}\left(\left( \lambda+\mu \, |(z,t)|^{\alpha q}\right)^{1/s}|f|^p\right)^{s} dz\,dt\right\}^{1/s} \left\{ \int_{\mathbb{H}^n}\left(\left( \lambda+\mu \, |(z,t)|^{\alpha q}\right)^{-1/s}\right)^{s'} dz\,dt\right\}^{1/s'}\right)^{s}\\
    &=\lVert \left( \lambda+\mu \, |(z,t)|^{\alpha q} \right)^{1/s}|f|^p\rVert^{s}_{s} \;\lVert \left( \lambda+\mu \, |(z,t)|^{\alpha q} \right)^{-1/s} \rVert^{s}_{s'}\\
    &=\left( \lambda \, \lVert f \rVert_q^q+\mu \, \lVert \,|(z,t)|^{\alpha}f \rVert_q^q\right)\;\lVert \left( \lambda+\mu \, |(z,t)|^{\alpha q} \right)^{-1/s} \rVert^{s}_{s'},
    \end{align*}
    where $s$ and $s'$ are dual H\"older exponents.
    
    Now, we calculate the second quantity on RHS using polar decomposition \eqref{polar_decomp_koranyi} on \(\mathbb{H}^n.\)
\begin{align*}
    \left(\int_{\mathbb{H}^n}\left( \lambda+\mu \, |(z,t)|^{\alpha q}\right)^{-1/(s-1)} dz \, dt \right)^{s-1}&=\left( \int_{0}^{\infty}\int_{S_K}\left( \lambda+\mu \, r^{\alpha q}\right)^{-1/(s-1)}r^{Q-1}\,d\sigma \,dr\right)^{s-1}\\
    &\asymp  \lambda^{-1}\left( \la/\mu \right)^{Q\beta/\alpha}.
\end{align*}
This completes the proof.
\end{proof}

The next lemma is instrumental in establishing the main theorem of this section.

\begin{lemma}\label{lem product inequality}
    Suppose that \(p,q,\alpha\) and \(\beta\) satisfies the hypothesis of Proposition \ref{prop sum inequality}. Then the following inequality
    $$\|f\|_p \lesssim \|f\|_q^{(1 - Q\beta/\alpha)} \; \| \, |(z,t)|^\alpha \, f\|_q^{Q\beta/\alpha}$$ holds for all $f \in L^q(\mathbb{H}^n).$
\end{lemma}

\begin{proof}
    Let $r > 0.$ If we consider the function $\delta_r(f),$ then it is easy to verify that
    $$\|\delta_r(f)\|_p = r^{-Q/p} \|f\|_p \;\; \text{and} \;\; \|\,|(z,t)|^\alpha \delta_r(f)\|_q  = r^{-(\alpha + Q/q)} \|\,|(z,t)|^\alpha f\|_q.$$
    By using Proposition \ref{prop sum inequality} for $\delta_r(f)$ in place of $f,$ we have
    \begin{equation}\label{value of r}
        \|f\|_p^q \lesssim \; \lambda^{-1} \; (\lambda/\mu)^{Q \beta/\alpha} \; r^{Q\beta q}\left( \lambda \, \|f\|_q^q + \mu \, r^{- \alpha q} \, \| \, |(z,t)|^\alpha f\|_q^q \right).
    \end{equation}
    Now, choose 
\[
r
= \left( \frac{\mu(\alpha - Q\beta)}{\lambda Q \beta} \right)^{1/(\alpha q)}
\, \|f\|_{q}^{-1/\alpha}
\, \big\| |(z,t)|^{\alpha} f \big\|_{q}^{1/\alpha}
>0
\]
in \eqref{value of r}. Substituting this choice into the preceding estimate and simplifying yields the desired inequality.

\end{proof}

We proceed to present the local uncertainty principle of Price for the Strichartz Fourier transform.

\begin{theorem}
    Let $E$ be a subset of $\Omega \times \mathbb{C}^n$ and $\alpha > Q/2.$ Then for all $f \in L^2(\mathbb{H}^n),$ we have
    \begin{equation}\label{mainineq}
        \int_E \left|\left(\mathcal{F}(f)\right)(a,w)\right|^2 d\nu_2(a)\,dw \lesssim |E| \; \|f\|_2^{2-Q/\alpha} \; \| \, |(z,t)|^\alpha f\|_2^{Q/\alpha}.
    \end{equation}
\end{theorem}

\begin{proof}
     If the RHS of \eqref{mainineq} is infinite, then the inequality holds trivially. Therefore, we may assume that the RHS of \eqref{mainineq} is finite. Let $E$ be a subset of $\Omega \times \mathbb{C}^n$ such that $0 < |E| < \infty.$ It follows that
    \begin{equation}\label{ineq}
        \int_E \left|\left(\mathcal{F}(f)\right)(a,w)\right|^2 d\nu_2(a)\,dw \leq |E| \; \left(\sup_{\Omega \times \mathbb{C}^n} \left|\left(\mathcal{F}(f)\right)(a,w)\right|\right)^2.
    \end{equation}
   Next, by choosing $p = 1$ and $q = 2$ in Lemma \ref{lem product inequality}, we have 
    \begin{equation}\label{inequality}
        \|f\|_1 \lesssim \|f\|_2^{(1 - Q/2\alpha)} \; \| \, |(z,t)|^\alpha \, f\|_2^{Q/2\alpha},
    \end{equation}
    for all $f \in L^2(\mathbb{H}^n)$ and $\alpha > Q/2.$ Further, for any $f \in L^1(\mathbb{H}^n),$ it follows from \eqref{L1 to L inf bdd} that
    \begin{equation}\label{fourier bound}
        \sup_{(a,w) \in \Omega \times \mathbb{C}^n} \left|\left(\mathcal{F}(f)\right)(a,w)\right| \lesssim \|f\|_1.
    \end{equation}
    Thus, the desired inequality \eqref{mainineq} is a direct consequence of \eqref{ineq}, \eqref{inequality} and \eqref{fourier bound}. 
\end{proof}

Furthermore, by adapting the argument used for the case \(0 < \alpha<n/2\) in \(\mathbb{R}^n\) and employing the polar decomposition associated with the Korányi norm, we obtain an analogous result for the Strichartz Fourier transform for the case $0 < \alpha < Q/2.$ For brevity, we omit the detailed proof, and the precise statement is as follows.

\begin{theorem}
  Let $E$ be a subset of $\Omega \times \mathbb{C}^n$ and $0<\alpha < Q/2.$ Then for all $f \in L^2(\mathbb{H}^n),$ we have  
  $$
  \int_E \left|\left(\mathcal{F}(f)\right)(a,w)\right|^2 d\nu_2(a)\,dw\lesssim |E|^{2\alpha/Q} \, \lVert \, |(z,t)|^{\alpha}f \rVert_2^{2}.
  $$
\end{theorem}

%-----------------------------------------------------------------------------

\section{Beurling's theorem}\label{sec 6}

In 2020, Elloumi et al. \cite{generalized-beurling} proved a generalized version of Beurling's theorem on $\mathbb{R}^n.$ The precise statement of the theorem is as follows.

\begin{theorem}\label{generalizedbeurling}
    Let $p,q \geq 1, N \geq 0$ and $f \in L^2(\mathbb{R}^n).$ Suppose that 
    $$\int_{\mathbb{R}^n} \int_{\mathbb{R}^n} \frac{|f(x)|^p \, |\widehat{f}(y)|^q}{(1+\|x\|+\|y\|)^N} \, e^{(2\max\{p,q\}
-1)\|x\|\|y\|} \, dx\,dy < \infty.$$
    Then the following holds:
    \begin{enumerate}
        \item If $N > n,$ then $f$ is of the form $P(x) \, e^{-a\|x\|^2},$ where $a > 0$ and $P$ is a polynomial of degree $d$ such that $d < (N-n)/2.$
        \item If $N \leq n,$ then $f$ vanishes almost everywhere on $\mathbb{R}^n.$
    \end{enumerate}
\end{theorem}

In this section, we establish a weak form of Beurling's theorem for the Strichartz Fourier transform, restricting our attention to radial functions on \(\mathbb{H}^n\). Recall that a function \(f\) on \(\mathbb{H}^n\) is radial if it is invariant under the natural action of \(U(n)\); equivalently, \(f(z,t)=f_0(|z|,t)\) for a unique function \(f_0\) on \(\mathbb{R}^+ \times \mathbb{R}\). For such radial functions, the Strichartz Fourier transform admits a markedly simpler representation; see \cite[(7.30)]{Strichartz_Fourier_thangavelu}. Explicitly, 
$$
\mathcal{F}(f)((\lambda,(2k+n)|\lambda|),w)=c_{n,k}\,R^{n-1}_{k}(-\lambda,f)\,e^{n-1}_{k,\lambda}(w,0),
$$
where \(R^{n-1}_{k}(\la,f)\) is the \(k^{th}\) Laguerre coefficient of $$f^{\lambda}(z)=\int_{\mathbb{R}}f(z,t)\,e^{i \lambda t}\,dt$$ and given by 
$$
R^{n-1}_{k}(\la,f)=c_{n,k}\int_{\mathbb{C}^n}f^{\lambda}(z)\,\varphi^{n-1}_{k,\lambda}(z)\,dz.
$$

We now present an analogue of Beurling's theorem for the Strichartz Fourier transform.

\begin{theorem}
    Let $f \in L^2(\mathbb{H}^n)$ be a radial function on $\mathbb{H}^n$ and $N \geq 0.$ If $$\int_{\mathbb{H}^n} \int_\Omega \frac{|f(z,t)|^2 \; |\mathcal{F}(f)((\la,(2k+n)|\la|),0)|^2}{(1+|t|+|\la|)^N} \, e^{3|t||\la|} \, |\la|^{-n} \, \frac{k!}{\Gamma(k+n)} \, d\nu_2(a) \, dz \, dt < \infty,$$
    where $a = (\la,(2k+n)|\la|),$ then the following holds:
    \begin{enumerate}
    
        \item If $N > 1,$ then $$f(z,t) = e^{-at^2} \left(\sum_{i=0}^d \psi_i(|z|) \, t^i\right)$$
    where $a > 0, d< (N-1)/2$ and $\psi_i(y) \in L^2\big([0,\infty),y^{2n-1}dy\big)$ for $0 \leq i \leq d.$
    \item If $N \leq 1,$ then $f = 0.$
    \end{enumerate}
\end{theorem}

\begin{proof}
    Since $f$ is a radial function on $\mathbb{H}^n,$ we have 
    \begin{align*}
      \mathcal{F}(f)((\la,(2k+n)\la),0) &= c_{n,k} \left(c_{n,k} \int_{\mathbb{C}^n} f^{-\la}(z) \, \varphi_{k,-\la}^{n-1}(z) \, dz\right) e_{k,\la}^{n-1}(0,0) \\
       &= c_{n,k} \int_{\mathbb{C}^n} f^{-\la}(z) \, \varphi_{k,\la}^{n-1}(z) \, dz\\
       &= c_{n,k} \int_0^\infty \sigma(S^{2n-1}) \, f^{-\la}(y)\, \varphi_{k,\la}^{n-1}(y) \, y^{2n-1} dy \\
       &= c_{n,k} \; \sigma(S^{2n-1}) \int_0^\infty f^{-\la}(y) \, L_k^{n-1}\left(\frac{1}{2}|\la|y^2\right) e^{-\frac{1}{4}|\la|y^2} y^{2n-1} dy\\
       &= c_{n,k} \; \sigma(S^{2n-1}) \, \left(\frac{|\la|}{2}\right)^{-n/2} \left(\frac{2k!}{\Gamma(k+n)}\right)^{-1/2}\\&\hspace{0.85cm}\int_0^\infty f^{-\la}(y) \left(\frac{|\la|}{2}\right)^{n/2} \left(\frac{2k!}{\Gamma(k+n)}\right)^{1/2}L_k^{n-1}\left(\frac{1}{2}|\la|y^2\right) e^{-\frac{1}{4}|\la|y^2} y^{2n-1} dy.
    \end{align*}
    
    Since from \cite[Lemma 2]{huang-liu}, the system $\{\phi_{k,r}(y):=r^n \, \phi_k(ry): k \in \mathbb{N}\cup\{0\}\}$ forms an orthonormal basis of the space $L^2\big([0,\infty),y^{2n-1}dy\big),$ where $r = \left( \frac{|\la|}{2}\right)^{1/2} > 0$ and 
    $$\phi_k(y) = \left( \frac{2k!}{\Gamma(k+n)} \right)^{1/2} e^{-y^2/2} \, L_k^{n-1}(y^2),$$
    it follows that $$\mathcal{F}(f)((\la,(2k+n)\la),0) = c_{n,k} \; \sigma(S^{2n-1}) \, \left(\frac{|\la|}{2}\right)^{-n/2} \left(\frac{2k!}{\Gamma(k+n)}\right)^{-1/2} \langle f^{-\la},\phi_{k,r} \rangle.$$ 
    Thus, we have
    \begin{align}
        \int_0^\infty |f^{-\la}(y)|^2 \, y^{2n-1} \, dy &= \sum_{k=0}^\infty |\langle f^{-\la},\phi_{k,r} \rangle|^2 \nonumber\\
        &= 2 \, \sigma(S^{2n-1})^{-2} \left( \frac{|\la|}{2} \right)^n \sum_{k=0}^\infty  (c_{n,k})^{-2} \left(\frac{k!}{\Gamma(k+n)}\right) |\mathcal{F}(f)((\la,(2k+n)\la),0)|^2. \label{flambda-equality}
    \end{align}
So, writing \eqref{flambda-equality} together with the radiality assumption on \(f\), we conclude that

    \begin{align*}
        &\int_{\mathbb{H}^n} \int_\Omega \frac{|f(z,t)|^2 \; |\mathcal{F}(f)((\la,(2k+n)|\la|),0)|^2}{(1+|t|+|\la|)^N} \, e^{3|t||\la|} \, |\la|^{-n} \, \frac{k!}{\Gamma(k+n)} \, d\nu_2(a) \, dz \, dt\\
        &= (2\pi)^{-2n-1} \int_{\mathbb{H}^n} \int_\mathbb{R} \frac{|f(z,t)|^2}{(1+|t|+|\la|)^N}  \left(\sum_{k=0}^\infty (c_{n,k})^{-2} \left(\frac{k!}{\Gamma(k+n)}\right) |\mathcal{F}(f)((\la,(2k+n)|\la|),0)|^2 \right) |\la|^{n}\\ &\hspace{14cm}e^{3|t||\la|} \, d\la \, dz \, dt\\
        &= (2\pi)^{-2n-1} \, 2^{n-1} \sigma(S^{2n-1})^2 \int_{\mathbb{H}^n} \int_\mathbb{R} \frac{|f(z,t)|^2}{(1+|t|+|\la|)^N} \int_0^\infty |f^{-\la}(y)|^2 \, y^{2n-1} \, dy \; e^{3|t||\la|} \, dz \, dt \, d\la\\
        &= (2\pi)^{-2n-1} \, 2^{n-1} \sigma(S^{2n-1})^3 \int_\mathbb{R} \int_\mathbb{R} \frac{1}{(1+|t|+|\la|)^N} \left(\int_0^\infty \int_0^\infty |f(x,t)|^2 \, |f^{-\la}(y)|^2 \, x^{2n-1} \, y^{2n-1} \, dx \, dy\right) \\ &\hspace{14.5cm}e^{3|t||\la|} \, dt \, d\la \\
        &= (2\pi)^{-2n-1} \, 2^{n-1} \sigma(S^{2n-1})^3 \int_\mathbb{R} \int_\mathbb{R} \frac{G(t,\la)}{(1+|t|+|\la|)^N}\; e^{3|t||\la|} \, dt \, d\la,
    \end{align*}
    where $$G(t,\la) := \int_0^\infty \int_0^\infty |f(x,t)|^2 \; |f^{-\la}(y)|^2 \, x^{2n-1} \, y^{2n-1} \, dx \, dy.$$

    By given hypothesis, we have
    $$(2\pi)^{-2n-1} \, 2^{n-1} \sigma(S^{2n-1})^3 \int_\mathbb{R} \int_\mathbb{R} \frac{G(t,\la)}{(1+|t|+|\la|)^N}\; e^{3|t||\la|} \, dt \, d\la < \infty.$$

    The remainder of the proof follows analogously to that of \cite[Theorem 2]{huang-liu}, utilizing Theorem \ref{generalizedbeurling} with $p = q = 2$ and $n=1.$ Therefore, we omit the details.
\end{proof}

%-----------------------------------------------------------------------------
 
\section*{Acknowledgement}
    The first and third author gratefully acknowledge the Indian Institute of Technology Delhi for providing the Institute Assistantship.

\section*{Data Availability} 
    Data sharing does not apply to this article as no datasets were generated or analyzed during the current study.

\section*{Competing Interests}
    The authors declare that they have no competing interests. 

\bibliographystyle{acm}
\bibliography{ref}

@article {jaming-integral_operators,
    AUTHOR = {Ghobber, Saifallah and Jaming, Philippe},
     TITLE = {Uncertainty principles for integral operators},
   JOURNAL = {Studia Math.},
  FJOURNAL = {Studia Mathematica},
    VOLUME = {220},
      YEAR = {2014},
    NUMBER = {3},
     PAGES = {197--220},
      ISSN = {0039-3223,1730-6337},
   MRCLASS = {42C20 (42A38)},
  MRNUMBER = {3173045},
MRREVIEWER = {V.\ V.\ Peller},
       DOI = {10.4064/sm220-3-1},
       URL = {https://doi.org/10.4064/sm220-3-1},
}

@article {beals-gaveau-greiner,
    AUTHOR = {Beals, Richard W. and Gaveau, Bernard and Greiner, Peter C.
              and Vauthier, Jacques},
     TITLE = {The {L}aguerre calculus on the {H}eisenberg group. {II}},
   JOURNAL = {Bull. Sci. Math. (2)},
  FJOURNAL = {Bulletin des Sciences Math\'ematiques. 2e S\'erie},
    VOLUME = {110},
      YEAR = {1986},
    NUMBER = {3},
     PAGES = {225--288},
      ISSN = {0007-4497},
   MRCLASS = {22E25 (22E30 32F20 47G05 58G15)},
  MRNUMBER = {877171},
MRREVIEWER = {Steven\ George\ Krantz},
}

@book {beals1988calculus,
    AUTHOR = {Beals, Richard W. and Greiner, Peter C.},
     TITLE = {Calculus on Heisenberg manifolds},
    SERIES = {Annals of Mathematics Studies},
    VOLUME = {119},
 PUBLISHER = {Princeton University Press, Princeton, NJ},
      YEAR = {1988},
}

@article {rottensteiner-ruzhansky-harmonic-oscillators,
    AUTHOR = {Rottensteiner, David and Ruzhansky, Michael},
     TITLE = {Harmonic and anharmonic oscillators on the {H}eisenberg group},
   JOURNAL = {J. Math. Phys.},
  FJOURNAL = {Journal of Mathematical Physics},
    VOLUME = {63},
      YEAR = {2022},
    NUMBER = {11},
     PAGES = {Paper No. 111509, 23},
      ISSN = {0022-2488,1089-7658},
   MRCLASS = {35P05 (43A80)},
  MRNUMBER = {4503714},
       DOI = {10.1063/5.0106068},
       URL = {https://doi.org/10.1063/5.0106068},
}

@phdthesis{rottensteiner-thesis,
  title={Time-frequency analysis on the Heisenberg group},
  author={Rottensteiner, David},
  year={2014},
  school={Imperial College London}
}

@book {fisher-ruzhansky-book,
    AUTHOR = {Fischer, Veronique and Ruzhansky, Michael},
     TITLE = {Quantization on nilpotent {L}ie groups},
    SERIES = {Progress in Mathematics},
    VOLUME = {314},
 PUBLISHER = {Birkh\"auser/Springer, Cham},
      YEAR = {2016},
     PAGES = {xiii+557},
      ISBN = {978-3-319-29557-2; 978-3-319-29558-9},
   MRCLASS = {22E25 (22E30 35R03 35S05 43A80 46L05)},
  MRNUMBER = {3469687},
MRREVIEWER = {Antoni\ Wawrzy\'nczyk},
       DOI = {10.1007/978-3-319-29558-9},
       URL = {https://doi.org/10.1007/978-3-319-29558-9},
}

@article {Donoho-stark_paper,
    AUTHOR = {Donoho, David L. and Stark, Philip B.},
     TITLE = {Uncertainty principles and signal recovery},
   JOURNAL = {SIAM J. Appl. Math.},
  FJOURNAL = {SIAM Journal on Applied Mathematics},
    VOLUME = {49},
      YEAR = {1989},
    NUMBER = {3},
     PAGES = {906--931},
      ISSN = {0036-1399},
   MRCLASS = {42A05 (94A11)},
  MRNUMBER = {997928},
MRREVIEWER = {A.\ Bultheel},
       DOI = {10.1137/0149053},
       URL = {https://doi.org/10.1137/0149053},
}

@article {Pritam-Lacunary_spherical_function,
    AUTHOR = {Ganguly, Pritam and Thangavelu, Sundaram},
     TITLE = {On the lacunary spherical maximal function on the {H}eisenberg
              group},
   JOURNAL = {J. Funct. Anal.},
  FJOURNAL = {Journal of Functional Analysis},
    VOLUME = {280},
      YEAR = {2021},
    NUMBER = {3},
     PAGES = {Paper No. 108832, 32},
      ISSN = {0022-1236,1096-0783},
   MRCLASS = {43A80 (22E25 22E30 42B15 42B25)},
  MRNUMBER = {4170795},
MRREVIEWER = {Wenjuan\ Li},
       DOI = {10.1016/j.jfa.2020.108832},
       URL = {https://doi.org/10.1016/j.jfa.2020.108832},
}

@book{thangaveluheisenberg,
  title={Harmonic analysis on the Heisenberg group},
  author={Thangavelu, Sundaram},
  volume={159},
  year={2012},
  publisher={Springer Science \& Business Media}
}

@incollection {Strichartz_Fourier_thangavelu,
    AUTHOR = {Thangavelu, Sundaram},
     TITLE = {A scalar-valued {F}ourier transform for the {H}eisenberg
              group},
 BOOKTITLE = {From classical analysis to analysis on fractals. {V}ol. 1. {A}
              tribute to {R}obert {S}trichartz},
    SERIES = {Appl. Numer. Harmon. Anal.},
     PAGES = {137--163},
 PUBLISHER = {Birkh\"auser/Springer, Cham},
      YEAR = {(2023)},
      ISBN = {978-3-031-37799-0; 978-3-031-37800-3},
   MRCLASS = {43A80 (33C45 35R03 42C05)},
  MRNUMBER = {4676386},
       DOI = {10.1007/978-3-031-37800-3\_7},
       URL = {https://doi.org/10.1007/978-3-031-37800-3_7},
}

@book {wong_weyl-transforms,
    AUTHOR = {Wong, M. W.},
     TITLE = {Weyl transforms},
    SERIES = {Universitext},
 PUBLISHER = {Springer-Verlag, New York},
      YEAR = {1998},
     PAGES = {viii+158},
      ISBN = {0-387-98414-3},
   MRCLASS = {47G30 (42B10 43A32 44A15 81S30)},
  MRNUMBER = {1639461},
MRREVIEWER = {Kh\'elifa\ Trim\`eche},
}

@book {thangavelu-uncertainty-book,
    AUTHOR = {Thangavelu, Sundaram},
     TITLE = {An introduction to the uncertainty principle},
    SERIES = {Progress in Mathematics},
    VOLUME = {217},
      NOTE = {Hardy's theorem on Lie groups,
              With a foreword by Gerald B.\ Folland},
 PUBLISHER = {Birkh\"auser Boston, Inc., Boston, MA},
      YEAR = {2004},
     PAGES = {xiv+174},
      ISBN = {0-8176-4330-3},
   MRCLASS = {43A80 (33C45 33C55 42C10 43-02 43A20)},
  MRNUMBER = {2008480},
MRREVIEWER = {Weixing\ Zheng},
       DOI = {10.1007/978-0-8176-8164-7},
       URL = {https://doi.org/10.1007/978-0-8176-8164-7},
}

@book {szego,
    AUTHOR = {Szeg\"o, Gabor},
     TITLE = {Orthogonal {P}olynomials},
    SERIES = {American Mathematical Society Colloquium Publications},
    VOLUME = {Vol. 23},
 PUBLISHER = {American Mathematical Society, New York},
      YEAR = {1939},
     PAGES = {ix+401},
   MRCLASS = {42.3X},
  MRNUMBER = {77},
MRREVIEWER = {J.\ A.\ Shohat},
}

@book {folland-real_analysis_book,
    AUTHOR = {Folland, Gerald B.},
     TITLE = {Real analysis},
    SERIES = {Pure and Applied Mathematics (New York)},
   EDITION = {Second},
      NOTE = {Modern techniques and their applications,
              A Wiley-Interscience Publication},
 PUBLISHER = {John Wiley \& Sons, Inc., New York},
      YEAR = {1999},
     PAGES = {xvi+386},
      ISBN = {0-471-31716-0},
   MRCLASS = {00A05 (26-01 28-01 46-01)},
  MRNUMBER = {1681462},
}

@book {folland-harmonic_phase_space,
  title={Harmonic Analysis in Phase Space, 122},
  author={Folland, Gerald B.},
  year={1989},
  publisher={Princeton University Press, Princeton, NJ}
}

@article {generalized-beurling,
    AUTHOR = {Elloumi, M. and Baklouti, A. and Azaouzi, S.},
     TITLE = {A generalized {B}eurling theorem for some {L}ie groups},
   JOURNAL = {Math. Notes},
  FJOURNAL = {Mathematical Notes},
    VOLUME = {107},
      YEAR = {2020},
    NUMBER = {1-2},
     PAGES = {42--53},
      ISSN = {0001-4346,1573-8876},
   MRCLASS = {43A30 (42B10 43A80)},
  MRNUMBER = {4144399},
       DOI = {10.1134/s0001434620010058},
       URL = {https://doi.org/10.1134/s0001434620010058},
}

@article {huang-liu,
    AUTHOR = {Huang, Jizheng and Liu, Heping},
     TITLE = {An analogue of {B}eurling's theorem for the {L}aguerre
              hypergroup},
   JOURNAL = {J. Math. Anal. Appl.},
  FJOURNAL = {Journal of Mathematical Analysis and Applications},
    VOLUME = {336},
      YEAR = {2007},
    NUMBER = {2},
     PAGES = {1406--1413},
      ISSN = {0022-247X,1096-0813},
   MRCLASS = {43A62},
  MRNUMBER = {2353023},
MRREVIEWER = {Sundaram\ Thangavelu},
       DOI = {10.1016/j.jmaa.2007.03.054},
       URL = {https://doi.org/10.1016/j.jmaa.2007.03.054},
}

@article {FSsurvey,
    AUTHOR = {Folland, Gerald B. and Sitaram, Alladi},
     TITLE = {The uncertainty principle: a mathematical survey},
   JOURNAL = {J. Fourier Anal. Appl.},
  FJOURNAL = {The Journal of Fourier Analysis and Applications},
    VOLUME = {3},
      YEAR = {1997},
    NUMBER = {3},
     PAGES = {207--238},
      ISSN = {1069-5869,1531-5851},
   MRCLASS = {42B10 (22E30 43A30 81S05 81S30)},
  MRNUMBER = {1448337},
MRREVIEWER = {Michael\ Cowling},
       DOI = {10.1007/BF02649110},
       URL = {https://doi.org/10.1007/BF02649110},
}

@article {Narayanan-Ratnakumar,
    AUTHOR = {Narayanan, E. K. and Ratnakumar, P. K.},
     TITLE = {Benedicks' theorem for the {H}eisenberg group},
   JOURNAL = {Proc. Amer. Math. Soc.},
  FJOURNAL = {Proceedings of the American Mathematical Society},
    VOLUME = {138},
      YEAR = {2010},
    NUMBER = {6},
     PAGES = {2135--2140},
      ISSN = {0002-9939,1088-6826},
   MRCLASS = {42B10 (22E30 43A05)},
  MRNUMBER = {2596052},
MRREVIEWER = {Franz\ Luef},
       DOI = {10.1090/S0002-9939-10-10272-X},
       URL = {https://doi.org/10.1090/S0002-9939-10-10272-X},
}

@article {arnal-ludwig,
    AUTHOR = {Arnal, Didier and Ludwig, Jean},
     TITLE = {Q.{U}.{P}.\ and {P}aley-{W}iener properties of unimodular,
              especially nilpotent, {L}ie groups},
   JOURNAL = {Proc. Amer. Math. Soc.},
  FJOURNAL = {Proceedings of the American Mathematical Society},
    VOLUME = {125},
      YEAR = {1997},
    NUMBER = {4},
     PAGES = {1071--1080},
      ISSN = {0002-9939,1088-6826},
   MRCLASS = {43A30 (22D20 22E27)},
  MRNUMBER = {1353372},
MRREVIEWER = {Jonathan\ M.\ Rosenberg},
       DOI = {10.1090/S0002-9939-97-03608-3},
       URL = {https://doi.org/10.1090/S0002-9939-97-03608-3},
}

@article {benedicks-original-paper,
    AUTHOR = {Benedicks, Michael},
     TITLE = {On {F}ourier transforms of functions supported on sets of
              finite {L}ebesgue measure},
   JOURNAL = {J. Math. Anal. Appl.},
  FJOURNAL = {Journal of Mathematical Analysis and Applications},
    VOLUME = {106},
      YEAR = {1985},
    NUMBER = {1},
     PAGES = {180--183},
      ISSN = {0022-247X},
   MRCLASS = {43A25 (42B10)},
  MRNUMBER = {780328},
MRREVIEWER = {S.\ Hartman},
       DOI = {10.1016/0022-247X(85)90140-4},
       URL = {https://doi.org/10.1016/0022-247X(85)90140-4},
}

@article {Sitaram-Sundari-Thangavelu-uncertainty_principles,
    AUTHOR = {Sitaram, A. and Sundari, M. and Thangavelu, S.},
     TITLE = {Uncertainty principles on certain {L}ie groups},
   JOURNAL = {Proc. Indian Acad. Sci. Math. Sci.},
  FJOURNAL = {Indian Academy of Sciences. Proceedings. Mathematical
              Sciences},
    VOLUME = {105},
      YEAR = {1995},
    NUMBER = {2},
     PAGES = {135--151},
      ISSN = {0253-4142,0973-7685},
   MRCLASS = {43A15 (22D15 22E25 43A30)},
  MRNUMBER = {1350473},
MRREVIEWER = {John\ F.\ Price},
       DOI = {10.1007/BF02880360},
       URL = {https://doi.org/10.1007/BF02880360},
}

@article {Price-Sitaram-transforms_on_LCT_groups,
    AUTHOR = {Price, John F. and Sitaram, Alladi},
     TITLE = {Functions and their {F}ourier transforms with supports of
              finite measure for certain locally compact groups},
   JOURNAL = {J. Funct. Anal.},
  FJOURNAL = {Journal of Functional Analysis},
    VOLUME = {79},
      YEAR = {1988},
    NUMBER = {1},
     PAGES = {166--182},
      ISSN = {0022-1236},
   MRCLASS = {43A30 (22D12 22E30 22E46)},
  MRNUMBER = {950089},
MRREVIEWER = {Takeshi\ Kawazoe},
       DOI = {10.1016/0022-1236(88)90035-3},
       URL = {https://doi.org/10.1016/0022-1236(88)90035-3},
}

@article {Vemuri-benedicks,
    AUTHOR = {Vemuri, M. K.},
     TITLE = {Benedicks' theorem for the {W}eyl transform},
   JOURNAL = {J. Math. Anal. Appl.},
  FJOURNAL = {Journal of Mathematical Analysis and Applications},
    VOLUME = {452},
      YEAR = {2017},
    NUMBER = {1},
     PAGES = {209--217},
      ISSN = {0022-247X,1096-0813},
   MRCLASS = {43A80 (22E30 42B10 81S30)},
  MRNUMBER = {3628015},
MRREVIEWER = {Ioannis\ R.\ Parissis},
       DOI = {10.1016/j.jmaa.2017.02.032},
       URL = {https://doi.org/10.1016/j.jmaa.2017.02.032},
}

@article {somnath-benedicks_motion_group,
    AUTHOR = {Ghosh, Somnath and Srivastava, R. K.},
     TITLE = {Benedicks-{A}mrein-{B}erthier theorem for the {H}eisenberg
              motion group},
   JOURNAL = {Bull. Lond. Math. Soc.},
  FJOURNAL = {Bulletin of the London Mathematical Society},
    VOLUME = {54},
      YEAR = {2022},
    NUMBER = {2},
     PAGES = {526--539},
      ISSN = {0024-6093,1469-2120},
   MRCLASS = {43A30 (43A80 44A35)},
  MRNUMBER = {4441627},
MRREVIEWER = {Rajab\ Ali\ Kamyabi Gol},
       DOI = {10.1112/blms.12582},
       URL = {https://doi.org/10.1112/blms.12582},
}

@article {jaming-Ghobber-bessel_transform,
    AUTHOR = {Ghobber, Saifallah and Jaming, Philippe},
     TITLE = {Strong annihilating pairs for the {F}ourier-{B}essel
              transform},
   JOURNAL = {J. Math. Anal. Appl.},
  FJOURNAL = {Journal of Mathematical Analysis and Applications},
    VOLUME = {377},
      YEAR = {2011},
    NUMBER = {2},
     PAGES = {501--515},
      ISSN = {0022-247X,1096-0813},
   MRCLASS = {44A15 (42A38 42B10 44A20)},
  MRNUMBER = {2769153},
MRREVIEWER = {H.\ S. P. Shrivastava},
       DOI = {10.1016/j.jmaa.2010.11.015},
       URL = {https://doi.org/10.1016/j.jmaa.2010.11.015},
}

@article {Negzaoui-Laguerre_transform_benedicks,
    AUTHOR = {Negzaoui, Selma and Rebhi, Sami},
     TITLE = {On the concentration of a function and its {L}aguerre-{B}essel
              transform},
   JOURNAL = {Math. Inequal. Appl.},
  FJOURNAL = {Mathematical Inequalities \& Applications},
    VOLUME = {22},
      YEAR = {2019},
    NUMBER = {3},
     PAGES = {825--836},
      ISSN = {1331-4343,1848-9966},
   MRCLASS = {42B10 (43A62)},
  MRNUMBER = {3986993},
MRREVIEWER = {Mohamed\ Sifi},
       DOI = {10.7153/mia-2019-22-57},
       URL = {https://doi.org/10.7153/mia-2019-22-57},
}

@article {Daher-Opdam_Cherednik_transform,
    AUTHOR = {Achak, Azzedine and Daher, R.},
     TITLE = {Benedicks-{A}mrein-{B}erthier type theorem related to
              {O}pdam-{C}herednik transform},
   JOURNAL = {J. Pseudo-Differ. Oper. Appl.},
  FJOURNAL = {Journal of Pseudo-Differential Operators and Applications},
    VOLUME = {9},
      YEAR = {2018},
    NUMBER = {2},
     PAGES = {431--441},
      ISSN = {1662-9981,1662-999X},
   MRCLASS = {42A38 (34B30 44A35)},
  MRNUMBER = {3804043},
MRREVIEWER = {Branko\ Sari\'c},
       DOI = {10.1007/s11868-017-0189-9},
       URL = {https://doi.org/10.1007/s11868-017-0189-9},
}

@article {Faris-inequalities_local_uuncertainty,
    AUTHOR = {Faris, William G.},
     TITLE = {Inequalities and uncertainty principles},
   JOURNAL = {J. Math. Phys.},
  FJOURNAL = {Journal of Mathematical Physics},
    VOLUME = {19},
      YEAR = {1978},
    NUMBER = {2},
     PAGES = {461--466},
      ISSN = {0022-2488,1089-7658},
   MRCLASS = {81.47},
  MRNUMBER = {484134},
MRREVIEWER = {Christiane\ Coudray},
       DOI = {10.1063/1.523667},
       URL = {https://doi.org/10.1063/1.523667},
}

@article {Price-inequalities_local_uuncertainty,
    AUTHOR = {Price, John F.},
     TITLE = {Inequalities and local uncertainty principles},
   JOURNAL = {J. Math. Phys.},
  FJOURNAL = {Journal of Mathematical Physics},
    VOLUME = {24},
      YEAR = {1983},
    NUMBER = {7},
     PAGES = {1711--1714},
      ISSN = {0022-2488,1089-7658},
   MRCLASS = {81C99 (81B05)},
  MRNUMBER = {709504},
MRREVIEWER = {Christiane\ Coudray},
       DOI = {10.1063/1.525916},
       URL = {https://doi.org/10.1063/1.525916},
}

@article {Price-sharp_local_uncertainty,
    AUTHOR = {Price, John F.},
     TITLE = {Sharp local uncertainty inequalities},
   JOURNAL = {Studia Math.},
  FJOURNAL = {Polska Akademia Nauk. Instytut Matematyczny. Studia
              Mathematica},
    VOLUME = {85},
      YEAR = {1986},
    NUMBER = {1},
     PAGES = {37--45},
      ISSN = {0039-3223,1730-6337},
   MRCLASS = {42B10 (46E35 81B05)},
  MRNUMBER = {879414},
MRREVIEWER = {Tatjana\ Ostrogorski},
       DOI = {10.4064/sm-85-1-37-45},
       URL = {https://doi.org/10.4064/sm-85-1-37-45},
}

@article {Smith-uncertainty_on_lct,
    AUTHOR = {Smith, Kennan T.},
     TITLE = {The uncertainty principle on groups},
   JOURNAL = {SIAM J. Appl. Math.},
  FJOURNAL = {SIAM Journal on Applied Mathematics},
    VOLUME = {50},
      YEAR = {1990},
    NUMBER = {3},
     PAGES = {876--882},
      ISSN = {0036-1399},
   MRCLASS = {94A11 (42A05 43A25)},
  MRNUMBER = {1050918},
MRREVIEWER = {Radomir\ S.\ Stankovi\'c},
       DOI = {10.1137/0150051},
       URL = {https://doi.org/10.1137/0150051},
}

@article {Tuan-Dohono_Star_for_Hankel,
    AUTHOR = {Tuan, Vu Kim},
     TITLE = {Uncertainty principles for the {H}ankel transform},
   JOURNAL = {Integral Transforms Spec. Funct.},
  FJOURNAL = {Integral Transforms and Special Functions. An International
              Journal},
    VOLUME = {18},
      YEAR = {2007},
    NUMBER = {5-6},
     PAGES = {369--381},
      ISSN = {1065-2469,1476-8291},
   MRCLASS = {46F12 (33C10 42A05 42A10 42A16 43A50 44A15)},
  MRNUMBER = {2326075},
MRREVIEWER = {Luis\ Daniel\ Abreu},
       DOI = {10.1080/10652460701320745},
       URL = {https://doi.org/10.1080/10652460701320745},
}

@article {Price-Sitaram-Local_uncertainty,
    AUTHOR = {Price, John F. and Sitaram, Alladi},
     TITLE = {Local uncertainty inequalities for locally compact groups},
   JOURNAL = {Trans. Amer. Math. Soc.},
  FJOURNAL = {Transactions of the American Mathematical Society},
    VOLUME = {308},
      YEAR = {1988},
    NUMBER = {1},
     PAGES = {105--114},
      ISSN = {0002-9947,1088-6850},
   MRCLASS = {22E30 (43A15)},
  MRNUMBER = {946433},
MRREVIEWER = {S.\ Sankaran},
       DOI = {10.2307/2000952},
       URL = {https://doi.org/10.2307/2000952},
}

@article {Atef-inequalities_on_laguerre_hypergroup,
    AUTHOR = {Atef, Rahmouni},
     TITLE = {Uncertainty inequalities on {L}aguerre hypergroup},
   JOURNAL = {Mediterr. J. Math.},
  FJOURNAL = {Mediterranean Journal of Mathematics},
    VOLUME = {10},
      YEAR = {2013},
    NUMBER = {1},
     PAGES = {333--351},
      ISSN = {1660-5446,1660-5454},
   MRCLASS = {43A62 (42B10)},
  MRNUMBER = {3019109},
MRREVIEWER = {Mohamed\ Sifi},
       DOI = {10.1007/s00009-012-0198-0},
       URL = {https://doi.org/10.1007/s00009-012-0198-0},
}

@article {Ciatti-Ricci-Sundari-uncertainity_on_stratified,
    AUTHOR = {Ciatti, Paolo and Ricci, Fulvio and Sundari, M.},
     TITLE = {Uncertainty inequalities on stratified nilpotent groups},
   JOURNAL = {Bull. Kerala Math. Assoc.},
  FJOURNAL = {Bulletin of Kerala Mathematics Association},
      YEAR = {2005},
     PAGES = {53--72},
      ISSN = {0973-2721},
   MRCLASS = {22E25},
  MRNUMBER = {2250035},
MRREVIEWER = {Sundaram\ Thangavelu},
}

@article {sarkarbeurling,
    AUTHOR = {Sarkar, Rudra P. and Sengupta, Jyoti},
     TITLE = {Beurling's theorem for {R}iemannian symmetric spaces. {II}},
   JOURNAL = {Proc. Amer. Math. Soc.},
  FJOURNAL = {Proceedings of the American Mathematical Society},
    VOLUME = {136},
      YEAR = {2008},
    NUMBER = {5},
     PAGES = {1841--1853},
      ISSN = {0002-9939,1088-6826},
   MRCLASS = {43A85 (22E30)},
  MRNUMBER = {2373616},
MRREVIEWER = {Dachun\ Yang},
       DOI = {10.1090/S0002-9939-07-08990-3},
       URL = {https://doi.org/10.1090/S0002-9939-07-08990-3},
}

@article {thangavelubeurling,
    AUTHOR = {Thangavelu, Sundaram},
     TITLE = {Beurling's theorem on the {H}eisenberg group},
   JOURNAL = {Ark. Mat.},
  FJOURNAL = {Arkiv f\"or Matematik},
    VOLUME = {60},
      YEAR = {2022},
    NUMBER = {2},
     PAGES = {417--442},
      ISSN = {0004-2080,1871-2487},
   MRCLASS = {43A85 (33C45 35P10 42C05)},
  MRNUMBER = {4500373},
MRREVIEWER = {Jinsen\ Xiao},
}

@article {Beurling-rudra-parui,
    AUTHOR = {Parui, Sanjay and Sarkar, Rudra P.},
     TITLE = {Beurling's theorem and {$L^p$}-{$L^q$} {M}organ's theorem for
              step two nilpotent {L}ie groups},
   JOURNAL = {Publ. Res. Inst. Math. Sci.},
  FJOURNAL = {Kyoto University. Research Institute for Mathematical
              Sciences. Publications},
    VOLUME = {44},
      YEAR = {2008},
    NUMBER = {4},
     PAGES = {1027--1056},
      ISSN = {0034-5318,1663-4926},
   MRCLASS = {22E30 (43A80)},
  MRNUMBER = {2477903},
MRREVIEWER = {Takeshi\ Kawazoe},
       DOI = {10.2977/prims/1231263778},
       URL = {https://doi.org/10.2977/prims/1231263778},
}

@article{Dasgupta-Gulia-Beurling_theorem,
  title={Some Versions of {B}eurling's Theorem on {H}-type Groups},
  author={Dasgupta, Aparajita and Gulia, Prerna and Pusti, Sanjoy and Thangavelu, Sundaram},
  journal={arXiv:2505.15186},
  year={2025}
}

@article {bonami,
    AUTHOR = {Bonami, Aline and Demange, Bruno and Jaming, Philippe},
     TITLE = {Hermite functions and uncertainty principles for the {F}ourier
              and the windowed {F}ourier transforms},
   JOURNAL = {Rev. Mat. Iberoamericana},
  FJOURNAL = {Revista Matem\'atica Iberoamericana},
    VOLUME = {19},
      YEAR = {2003},
    NUMBER = {1},
     PAGES = {23--55},
      ISSN = {0213-2230},
   MRCLASS = {42B10 (33C45 33C47 94A12)},
  MRNUMBER = {1993414},
MRREVIEWER = {H.\ S. P. Shrivastava},
       DOI = {10.4171/RMI/337},
       URL = {https://doi.org/10.4171/RMI/337},
}

\end{document}